\documentclass[11pt]{amsart}
\usepackage{aliascnt, quiver, amssymb}
\usepackage{hank} 
\usepackage[style=alphabetic]{biblatex}
\title{Rigidity for Families of Abelian Varieties over $\mathcal{M}_{g,n}$}
\author{Nathan H. Morris}

\renewcommand{\bar}[1]{\overline{#1}}
\renewcommand{\tilde}[1]{\widetilde{#1}}

\begin{document}
	
	\maketitle
	\begin{abstract}
		We provide a new proof of Farb's theorem on the global rigidity of the period map $\mathcal{M}_{g,n} \to \mathcal{A}_g$. We extend this classification to essentially include the period map for $\mathcal{M}_{2,n}$. We also provide a classification of holomorphic maps $\mathcal{M}_g \to \mathcal{A}_{h}$ for $h \leq 2g-1$. The only nonconstant maps in this range are given by taking the Jacobian of a curve and attaching a fixed Abelian variety. In addition to the group theory of the mapping class group, we make extensive use of the theory of variations of Hodge structures and the birational geometry on $\overline{\mathcal{M}}_2$. 
	\end{abstract}

\section{Introduction}

	\subsection{Background}
		Let $\shm_{g,n}$ and $\sha_g$ denote the moduli spaces of genus $g$ curves with $n$ marked points and $g$-dimensional Abelian varieties, respectively. We take as our starting point the following theorem of Farb:
		\begin{theorem}\cite[Theorem 1.1]{farbrigid}\label{farbrigid}
			Let $f: \shm_{g,n} \to \sha_h$ be a holomorphic map with $h \leq g$. Then $h=g$ and $f$ is the Torelli map. 
		\end{theorem}
		
		Thus, the period map is the unique way to holomorphically associate a $g$-dimensional abelian variety to a genus $g$ curve. Farb's proof relies on a group-theoretic classification of representations of the mapping class group due to Korkmaz:
		
		\begin{theorem}\cite[Theorem 2]{korkmaz}
			Let $\rho: \Mod(g,n) \to \GL(2g,\CC)$ be a representation of the mapping class group with $g \geq 3$. Then either $\rho$ is trivial or $\rho$ is conjugate to the symplectic representation $\Mod(g,n) \to \Sp(2g,\ZZ)$. 
		\end{theorem}
		
		The impetus of this paper was the desire for a more algebro-geometric, specifically Hodge-theoretic, proof of Farb's theorem. After all, it is ultimately a result about a weight one variation of Hodge structures. Specifically, since Hain provided a Hodge theoretic proof of \Cref{farbrigid} assuming \Cite[Theorem 2]{korkmaz} in \Cite[Section 3]{farbrigid}, the desire was to give a proof of Korkmaz's classification rooted in Hodge theory. 
		
		By exploiting the extension properties of period maps across the boundary of the Deligne-Mumford compactification of $\shm_{g,n}$, this was possible. However, in making the argument precise and having it hold for a map on $\shm_{g,n}$ as opposed to a finite cover, we still require a great deal of group theoretic input. Combining the rigid group theory of the mapping class group with the Hodge-theoretic package coming from a period map, we are able to repove the theorems of Farb and Korkmaz, extending the domain of applicability to include $\shm_{2,n}$, and the range to include Abelian varieties up to dimension $2g-1$. 
	\subsection{Statement of Results}
		Our main theorem is the following classification of low-dimensional families of Abelian varieties over $\shm_{g,n}$:
		\begin{theorem}\label{mga2g-1}
			Let $f: \shm_{g,n} \to \sha_h$ be a holomorphic map with $h \leq 2g-1$ and $g \geq 2$. If $h < g$ then $f$ is constant. When $h=g=2$, there are two distinct nonconstant maps, which, on points both associate to a genus 2 curve its Jacobian. If $h=g \geq 3$ then $f$ is the Torelli map. Otherwise, $f$ is either constant or is the Torelli map $\shm_{g,n} \to \sha_g$ followed by an inclusion to $\sha_h$ by attaching a constant factor. 
		\end{theorem}
		This theorem is sharp, in the sense that when $h=2g$ there is a new map given by the two copies of the Torelli map. It is also worth nothing that, at least when $g =2$, the fact that the representation underlies a weight one variation of Hodge structure is crucial. Indeed, there is a nontrivial 6-dimensional representation of $\Mod(2,n)$ given by the exterior square of the representation on homology. It has an invariant line, and the quotient gives a nontrivial representation $\Mod(2,n) \to \SO(3,2)$. We will say more about the two maps when $g=2$ below.
		
		We prove our theorem by induction on $g$. Our induction starts at $g=2$, for which the Torelli map had not before been shown to be rigid:
		\begin{theorem}\label{m2a2}
			Let $f: \shm_{2,n} \to \sha_2$ be a nonconstant holomorphic map with unipotent monodromy around the divisor $\Delta_0$. Then $f$ is the Torelli map. 
		\end{theorem}
		
		We give two proofs of this theorem. The first follows Farb's strategy, and to execute his program we give a classification of rank 4 representations of $\Mod(2,n)$:
		\begin{theorem}\label{mgagreps}
			Let $\rho: \Mod(g,n) \to \GL(2g,\CC)$ have infinite image. If $g=2$, $\rho$ is conjugate to $\std \otimes \chi$, where $\std: \Mod(2,n) \to \Sp(4,\ZZ)$ is the homology representation and $\chi: \Mod(2,n) \to \ZZ/10\ZZ$ is a character. If $g \geq 3$, then $\rho$ is conjugate to $\std: \Mod(g,n) \to \Sp(2g,\ZZ)$.
		\end{theorem}
		
		The second proof is more geometric, and essentially comes down to the fact that $\shm_2$ is affine, while $\shm_{2}^{ct} \simeq \sha_2$; the fact that $\sha_2$ contains complete curves means a map $\shm_2 \to \sha_2$ must extend across the boundary of the Deligne-Mumford compactification. The process of showing a map $\shm_{2,n} \to \sha_2$ factors through the projection is more subtle, and borrows many techniques from the group-theoretic work done earlier in the paper.
		
		The base case of the induction also requires the classification of maps $\shm_{2,n} \to \sha_3$. Again, the group theoretic work done earlier in the paper eliminates many possible cases of how such a map can degenerate across the boundary of $\shm_{2,n}$, but not all of them. We analyze the remaining cases using some of the birational geometry of $\bar{\shm_2}$, whose nef cone is explicitly known. By pulling back the ample class on $\sha_3$ to the boundary strata we rule out the remaining cases. 
		
		The presence of finite abelianization in the genus 2 mapping class group explains the nontrivial character twists of the standard representation. Because there are then two homotopy classes of maps $\shm_{2} \to \sha_2$, the Torelli map is not rigid as a map of stacks:
		\begin{theorem}\label{nostacky}
			Let $f: \shm_{2,n} \to \sha_2$ be a map of Deligne-Mumford stacks. Then $f$ is either the Torelli map, or a twist of the Torelli map by $-I$. 
		\end{theorem}
		These maps are the same at the level of points, hence induce the same map on coarse moduli spaces, but they act differently on automorphisms. 
		
		En route to proving to proving the classification of rank-2g representations of the mapping class group, we recount the proof of \Cite[Theorem 1]{korkmaz}:
		
		\begin{theorem}\label{nolowrankreps}
			Let $\rho: \Mod(g,n) \to \GL(r,\CC)$ with $r \leq 2g-1$. If $g = 2$ then the image of $\rho$ is a quotient of $\ZZ/10\ZZ$, and if $ \geq 3$ the $\rho$ is trivial. 
		\end{theorem}	
		
		We prove this with the same strategy as Korkmaz, by considering all the possible eigenvalues and Jordan forms of such a represenation, and showing that they all annihilate the commutator subgroup of $\Mod(g,n)$. By using some techniques from Hodge theory and an inductive strategy, we are able to simplify the proof somewhat. 
	
	\subsection{Motivation and History}
		Rigidity has been a subject of interest in Hodge theory for a number of decades. As Schmid notes in \Cite{schmid}, various rigidity theorems families of Abelian varieties, and later general variations of Hodge structure, were considered by Grothendieck, Arakelov, and Griffiths. Schmid's formulation of the rigidity theorem (\Cite[7.24]{schmid}), is the criterion Hain uses to prove Farb's theorem, and Farb exploits the same boundedness property of period domains in his proof. 
		
		\subsubsection{Local Forms of Rigidity} 
		Most treatments of rigidity however, have been more local in nature. One of the earlier systematic treatments of rigidity for variations of Hodge structures came from Peters \cite{petersrigid}. He considered period maps as holomorphic maps $X \to D/\Gamma$, and he identified the tangent space to the space of such maps, and within that morphism space, the subspace of deformations which remain Griffiths-transverse. For a Shimura variety like $\sha_g$, every tangent direction is Griffiths transverse, but this is not the case for most period maps. In any case, Peters' form of rigidity is a kind of \textit{infinitesimal} rigidity, and, for example, his work shows that any variation of Hodge structures whose algebraic monodromy group is simple (one realization of "big monodromy") is infinitesimally rigid as a map of varieties.
		
		In another direction, one can relate rigidity of a variation of Hodge structures to the rigidity of its underlying representation:
		\begin{conjecture}\cite[Conjecture 6.12]{littmotives}
			Let $g \geq 3$. Any irreducible representation $\rho$ of $\Mod(g,n)$ satisfies $H^1(\Mod(g,n), \ad \rho) = 0$
		\end{conjecture}
		
		Since $H^1(\ad \rho)$ is the tangent space at $\rho$ to the character variety of representations of $\Mod(g,n)$, this says that $[\rho]$ is a reduced isolated point. This form is stronger than Peters' rigidity. As Litt notes, the infinitesimal rigidity of representations of $\Mod(g,n)$ is closely related to Ivanov's conjecture about the abelianization of finite-index subgroups of the mapping class group. 
		
		It is interesting to remark that Farb and Hain's proofs of global rigidity do not say anything about infinitesimal rigidity of the underlying representation. That is, the Torelli map, as a point of $\Hom(\Mod(g,n), \GL(2g,\CC))/\GL(2g)$ can be the unique point of the character variety whether or not that point is reduced. However, by work of Hain, the symplectic representation is known to be infinitesimally rigid, and we make use of this in our proof, as it allows us to assume any representation of the mapping class group in our regime supports a complex variation of Hodge structure; when that representation is then shown to be rigid in the character variety, we'll know that the representation we started with was a VHS to begin with. 
		
		\subsubsection{Global Rigidity and Representation Classification}
		For years, the state of the art for classification of low dimensional representations was due to Korkmaz \Cite{korkmaz}, who classified mapping class group representations of rank up to 2g, and Franks-Handel who independently proved the triviality of representations of rank up to 2g-1 \Cite{FranksHandel2013}. Farb developed an argument for enhancing this classification of representations to a classification of morphisms $\shm_{g,n} \to \sha_g$. Later, Serv\'an proved a similar global rigidity statement for the Prym morphism on a cover of $\shm_g$ \cite{servanprym}. 
		
		The clearest antecedent to this paper is the work of Kaufmann-Salter-Zhang-Zhong \Cite{kaufmann2026linearrepresentationsmappingclass}, which provides a classification of representations of the mapping class group $\Mod(g,n) \to \GL(r,\CC)$ with $n \in \{0,1\}$, $g \geq 4$, and $r \leq 3g-3$. Interestingly, they note a nontrivial representation of $\Mod(g,1)$ of rank $2g+1$, given by the action on the homology of the unit tangent bundle. This representation does not occur in our classification because it is not semisimple, hence it cannot support a pure variation of Hodge structures. This representation, however, and its generalizations when $n>1$, is closely related to normal functions and maps from $\shm_{g,n}$ to the universal abelian variety $\mathcal{X}_g$, see \cite{xiyanrigid}. As they prove their result on induction on $g$, starting with $g=3$ and $r=2g$ as a base case, we expect \Cref{mgagreps} with some additional argument to slightly extend the rank range for which their result is valid. 
		
		While this paper was in preparation, Zhong \cite{xiyanrigid} proved the rigidity of the Torelli map on the moduli space of triple covers of genus $g \geq 4$ curves, and also provided a classification of all bi-affine representations of $\Mod(g,n)$ of rank at least $2g$. These latter representations are certain extensions of the symplectic representation, and have a natural interpretation in terms of mixed Hodge structures. Using this classification, Zhong, joint with Z. Zhang, announced a classification of maps $\shm_{g,n} \to \sha_h$ whose underlying representation has rank at most $3g-3$, using different, more topological, methods. 
		\subsubsection{Questions}
		Many conjectures about mapping class groups are motivated by the analogy between mapping class groups and lattices in Lie groups. Under this analogy, the standard representation on the homology of the surface is a sort of fundamental representation, in that our regime, it and the trivial representation constitute all low dimensional (symplectic) representations of the mapping class group. The next nontrivial representations we know of, in order of rank, include tensor powers of the standard representation, whose ranks are polynomial in $g$, and then the action on $H_1(\Sigma_{g,n}, \ZZ/2\ZZ)$, whose rank is exponential in $g$. It is natural to ask:
		
		\begin{question}
			Let $g \geq 2$. What is the lowest rank $r$ for which there is a nontrivial representation $\Mod(g,n) \to \GL(r,\CC)$ that is not constructed functorially from the standard representation?
		\end{question}
		
		The action of the mapping class group on homology of the surface with Abelian coefficients is closely related to covering constructions. These representations clearly preserve a symplectic form. Curt McMullen asked if all period maps $\shm_g \to \sha_h$ come from covering constructions. Given the results of Kaufmann-Salter-Zhang-Zhong, one cannot literally pose the same question for arbitrary representations $\Mod(g,n) \to \GL(r,\CC)$. However, the question may be well-posed if one only considers semisimple representations:
		\begin{question}
			Let $g \geq 2$ and let $\rho: \Mod(g,n) \to \GL(r,\CC)$ be a representation and let $\rho^{ss}$ be its semisimplification. Does $\rho^{ss}$ come from covering constructions?
		\end{question}
		
		We also have a question of methods. If, as we expect, the next nontrivial variation of Hodge structure on $\shm_{g,n}$ is given by the action on $H_1(\Sigma_{g,n}, \ZZ/2\ZZ)$, how might one prove it? Our methods rely on a linear induction, which is clearly inadequate if the next goalpost has rank exponential in $g$. So, given the unipotent part of the local monodromy $I+N$ of a variation of Hodge structures around $\Delta_0$, one needs a way to rule out most large swaths of nilpotent ranks. 
	\subsection{Outline}
		Section 2 gives some background on Hodge theory and variations of Hodge structures, and Section 3 gives background on the mapping class group. Experts can probably skip these sections. Section 4 is the group-theoretic core of the paper. It contains proofs of \Cref{nolowrankreps} and \Cref{mgagreps}. Section 5 recounts how Farb and Hain use a classificaiton of representations to prove global rigidity of the period map, proving \Cref{m2a2} and \Cref{nostacky}. Section 6 contains the second proof of \Cref{m2a2}, and proves \Cref{mga2g-1}. 
	\subsection{Acknowledgements}
		I thank Ben Bakker for many useful comments and discussions related to this article.  
	\subsection{A.I. Disclosure} 
		The author used ChatGPT free, and later plus, to do matrix computations while formulating arguments.
\section{Hodge-Theoretic Preliminaries}
	Let $k = \ZZ, \QQ$. A pure $k$-Hodge structure of weight $n$ is a finitely generated $k$-module $V$ such that $V_\CC := V \otimes_k \CC$ has a decreasing filtration $F^\bullet$ such that $$V_\CC= \bigoplus_{p+q=n} F^p V \cap \overline{F^q V}.$$ 
	
	We make free use of the Riemann-Hilbert correspondence, establishing the equivalence of categories between representations of the fundamental group of a variety, and holomorphic vector bundles with flat connections. 
	
	A polarizable complex variation of Hodge structures of weight $k$ on a smooth irreducible complex variety $X$ is a flat vector bundle $(E, \nabla)$ with a decreasing decreasing filtrations $F^\bullet$ and $F^{'\bullet}$, the latter anti-holomorphic, each satisfying Griffiths transversality: $$\nabla(F^p E) \subset F^{p-1} E \otimes \Omega_X^1$$
	$$\overline{\nabla}(F^{'p}E) \subset F^{'p-1} \otimes \overline{\Omega}_X^1.$$
	 Moreover, there exists a flat hermitian form $\psi$ on $E$ such that $\gr_F^p \gr_{F'}^q E = 0$ for $p + q \neq k$, and $(-1)^p \psi$ is positive definite on $\gr_F^p E$. 
	
	One deep result about variations of Hodge structure we will invoke many times is Deligne's semi-simplicity theorem:
	\begin{proposition}\cite[1.11-12]{Del87}\label{semisimplicity}
		A local system $\VV$ which underlies a polarizable complex variation of Hodge structure is semisimple. Moreover, $\VV$ may be canonically decomposed as
		$$\VV \simeq \bigoplus_i \UU_i \otimes W_i,$$ where the $\UU_i$ are pairwise non-isomorphic irreducible complex local systems which underlie a polarizable complex variation of Hodge structures, and the $W_i$ are complex vector spaces. The structures on each $\UU_i$ and $W_i$ are both unique up to shifts in the grading. 
	\end{proposition}
	
	From this, we obtain the following corollary which will apply to our rigidity theorems:
	\begin{corollary}\label{Delignerigidity}
		Let $\VV$ be a local system whose algebraic monodromy group is simple and acts irreducible on $\VV$. Then $\VV$ underlies at most one polarizable complex variation of Hodge structures, unique up to shifts in the grading. 
	\end{corollary}

	Delgine's proof of \Cref{semisimplicity} makes extensive use of the theorem of the fixed part, which we will also repeatedly invoke:
	\begin{proposition}\cite[Theorem 11]{petersvhs}
		Let $\VV$ be a $\QQ$-variation of Hodge structures of weight $k$ on an irreducible compactifiable complex variety $X$. Then $H^0(X, \VV) = \VV^{\pi_1(X)}$ has a weight $k$ Hodge structure, and for any $x \in X$, the map $H^0(X, \VV) \to \VV_x$ is a morphism of Hodge structures. 
	\end{proposition}

	\subsection{Limit Mixed Hodge Structures}
	The basic constructions of the limiting Hodge structure can be found in \cite{schmid}. Suppose $(E,\nabla)$ is a variation of Hodge structures on a smoooth complex variety $X$ which is the complement of a simple normal crossings divisor. Let $\rho: \pi_1(X,x) \to \GL(r,\CC)$ denote the corresponding representation. Fix an irreducible component $D$ of the divisor and suppose locally that $D$ is given by the equation $s = 0$. Let $T$ be the monodromy around $D$. Then $T$ is known to be quasi-unipotent: $T = T_{s} T_u$ with $T_u = I+N$ \cite[(4.5) Lemma]{schmid}. One may use $N$ to construct a filtration $W_\bullet$ on $\CC^r$, called the \textit{weight filtration}. It is unique, each quotient $\gr^W_k$ is a pure weight $k$ Hodge structure, and $N^r$ defines an isomorphism $\gr^W_{k+r} E \to \gr^W_{k-r}$. If the weight of the Hodge structure on $\rho$ is $k$, then $N^{k+1} = 0$. In this context, if $k+1$ is the minimal exponent which annihilates $N$, we say that $N$ has order of nilpotency $k+1$. 
	
	We now consider the local fundamental group around $D$. Let $N^\times$ be a punctured tubular neighborhood around $D$. Then $N^\times$ has the homotopy type of an $S^1$-fibration over $D$. The inclusion $\pi_1(N^\times) \to \pi_1(X)$ followed by $\rho$ gives lets one analyze the degeneration of the variation of Hodge structures near $D$.
	
	\begin{lemma}
		In the notation above, the image of $\rho(\pi_1(N^\times))$ is contained in the parabolic subgroup $P \subset \GL(r,\CC)$ of transformations preserving $W_\bullet$.
	\end{lemma}
		\begin{proof}
			Any element of $\pi_1(N^\times)$ commutes with the element $\gamma$ giving the local monodromy around $D$, hence $\rho(\pi_1(N^\times))$ preserves the nilpotent part $N$ of $T=\rho(\gamma)$. The weight filtration is constructed from $N$ and is unique. In particular, $W_\bullet$ is preserved by conjugation, hence it is preserved by elements commuting with $N$. 
		\end{proof}
	
	If the monodromy around $D$ is unipotent, then $\rho(T)$ lands in the unipotent radical of $P$, and there is an induced variation of Hodge structures $\psi: \pi_1(D) \to L$, where $L \simeq P/U$, given by the action of $\pi_1(D)$ on the vector space graded by the weight filtration. We have the following diagram: 
	\[\begin{tikzcd}
		1 & {\langle T \rangle} & {\pi_1(N^\times)} & {\pi_1(D)} & 1 \\
		1 & U & P & L & 1
		\arrow[from=1-1, to=1-2]
		\arrow[from=1-2, to=1-3]
		\arrow[from=1-2, to=2-2]
		\arrow[from=1-3, to=1-4]
		\arrow["\rho", from=1-3, to=2-3]
		\arrow[from=1-4, to=1-5]
		\arrow["\psi", from=1-4, to=2-4]
		\arrow[from=2-1, to=2-2]
		\arrow[from=2-2, to=2-3]
		\arrow[from=2-3, to=2-4]
		\arrow[from=2-4, to=2-5]
	\end{tikzcd}\]
	
	We will use this construction in the case where $X = \shm_{g,n}$ and $D = \Delta_0 \simeq \shm_{g-1,2}/S_2$. Here, the monodromy around $D$ is given by the Dehn twist about a nonseparating simple closed curve (see Section 2 below). We can then inductively extend a classification of representations of $\Mod(g-1,n+2)$ to a classification of representations on $\Mod(g, n)$. 

	\subsection{Deformations of Representations}
	We will use some results about deformations of representations in our argument. Most importantly, we have the following result of Mochizuki, which states that a representation can always be deformed to one underlying a polarizable complex variation of Hodge structure.
	\begin{theorem}\cite[Theorem 10.5]{Moc}\label{MocDeform}
		Let $\bar{X}$ be a smooth projective variety and let $D \subset \bar{X}$ be a strict normal crossings divisor, and let $X = \bar{X} \setminus D$. Any representation $$\rho:\pi_1(X) \to \GL(r,\CC)$$ can be deformed with constant determinant to a representation underlying a polarizable complex variation of Hodge structure. 
	\end{theorem}
	
	Additionally, we will use the following result of Bakker-Brunebarbe-Tsimerman which will allow us to deform a representation $\rho$ to one which underlies a polarizable complex variation of Hodge structures, which additionally has quasi-unipotent monodromy at infinity.
	
	\begin{theorem}\cite[Immediate Consequence of Corollary 8.24]{bakker2024linearshafarevichconjecturequasiprojective}\label{qudensity}
		Every component of the moduli stack $\shm_B(X,r)$ parameterizing conjugacy classes of rank $r$ reprsentaitons of $\pi_1(X)$ contains a point parameterizing a local system which underlies a polarizable complex variation of Hodge structures with quasi-unipotent monodromy at infinity. 
	\end{theorem}
	
	In particular, any isolated point of $\shm_B(X,r)$ underlies a polarizable complex variation of Hodge structures with quasi-unipotent monodromy around $D$. In our proofs, we will use these results to reduce to the case when the representation underlies a polarizable complex variation of Hodge structure with quasi-unipotent monodromy, classify the resulting representation, and then show the representation we found admits no nontrivial deformations, so that our initial deformation was in fact trivial. To make this precise we use the language of cohomological rigidity. 
	
	\begin{definition}
		Let $\bar{X}$ be a smooth projective variety and let $D \subset \bar{X}$ be a strict normal crossings divisor.  Let $X = \bar{X}\setminus D$. Let $$\rho: \pi_1(X) \to \GL(r,\CC)$$ be a representation with quasi-unipotent monodromy around $D$. $\rho$ is cohomologically rigid if $$H^1(\bar{X} j_{!*} \ad \rho) = 0,$$ where $j: X \injto \bar{X}$ is the inclusion, and $\rho$ is strongly cohomologically rigid if $$H^1(\pi_1(X), \ad \rho) = 0.$$
	\end{definition}
	
	From the Leray spectral sequence, one sees that if $\rho$ is strongly cohomologically rigid, then it is cohomologically rigid. Geometrically, cohomological rigidity means that a representation is a reduced isolated point in its character variety. 
	
	We record the following cohomologically rigid representations:
	\begin{lemma}\label{rigidadrho}
		\begin{enumerate}
			\item Let $g \geq 3$. Let $\rho: \Mod(g,n) \to \Sp(2g,\ZZ)$ be the standard homology representation. Then $\rho$ is cohomologically rigid. 
			\item Let $g \geq 2$. Then if $\rho: \Mod(g,n) \to \GL(r, \CC)$ is a representation such that the commutator subgroup of $\Mod(g,n)$ is in the kernel of $\rho$, then $\rho$ is cohomologically rigid. 
		\end{enumerate}
	\end{lemma}
		\begin{proof}
			For both items, we show strong cohomological rigidity. Item (1) follows from \cite[Lemma 10.3]{Hain}. For (2), we consider the inflation-restriction sequence associated to $G = \Mod(g,n)$ and $N = [G,G]$. Let $Q =G/N$ Note that when $g = 2, Q \simeq \ZZ/10\ZZ$ while when $g \geq 3$, $Q = 1$. Let $V = \ad \rho$. Then we have an exact sequence
			$$0 \to H^1(Q, V^N) \to H^1(G,V) \to H^1(N, V)^{Q} \to H^2(Q, V^N).$$
			For $g \geq 2$, the commutator subgroup $N$ is perfect, so $H^1(N,V) = \Hom(N, V) = 0$. If $g \geq 3,$ it is clear that $H^1(Q, V^N) = 0$. When $g = 2$, this is $H^1(\ZZ/10\ZZ, V^N)$. Since $N$ acts trivially on $\ad \rho$, $V^N = V$, so this is just $H^1(\ZZ/10\ZZ, V)$. By Maschke's theorem, higher cohomology of finite groups with coefficients in a vector space are trivial. 
		\end{proof}
	
	Hence, if any variation of Hodge structures on $\shm_{g,n}$ is conjugate to the symplectic representation or annihilates the commutator subgroup, then this is in fact true for any representation of $\Mod(g,n)$ for that rank. 
\section{Group Theoretic Preliminaries}
	Let $S = \Sigma_{g,n}$ denote a genus $g$ surface with $n$ marked points. We let $\Mod(g,n)$ denote the \textit{mapping class group} of orientation-preserving diffeomorphisms of $S$ fixing the marked points, up to isotopy. By a \textit{simple closed curve} on $S$ we mean an embedding $S^1 \to S$. A simple closed curve $c$ is \textit{separating} if $S \setminus c$ is disconnected, and is \textit{nonseparating} otherwise. 
	
	Let $c$ be a simple closed curve in $S$. A tubular neighborhood of $c$ is homeomorphic to the cylinder $S^1 \times I$. An orientation preserving diffeomorphism of $S$ which, on this cylinder, is of the form $(s,t) \mapsto (se^{2\pi t},t)$ is called a \textit{Dehn twist} about $c$, denoted $T_c$. By the Dehn-Lickorish theorem \cite{farbprimer}[Theorem 4.1], $\Mod(g,n)$ is generated by finitely many Dehn twists about nonseparating simple closed curves. 
	
	\subsection{Generators and Relations in the Mapping Class Group}
	While not strictly necessary for our proofs, it will be nice to fix an explicit generating set for $\Mod(g,n)$ for the remainder of the paper. We will work with the Lickorish generators. Let $a_1,b_1, \dots, a_g,b_g$ denote a standard homology basis of simple closed curves on a closed genus $g$ surface, where $a_i \cap b_i = \{e\}$ and otherwise the curves are disjoint. Let $c_1, \dots, c_{g-1}$ denote nonseparating simple closed curves with $c_i \cap b_i = c_i \cap b_{i+1} = \{e\}$ and otherwise disjoint. Finally, let $a_1 = e_0, e_1, \dots, e_{n}$ denote nonseparating simple closed curves homologous to $a_1$, such that the $i$th marked point lies in the collar bounded by $e_i$ and $e_{i+1}$. 
	
	\begin{theorem}\cite[Theorem 4.13 + $\varepsilon$]{farbprimer}
		$\Mod(g,n)$ is generated by the Dehn twists about the isotopy classes of $a_1, \dots, a_g, b_1, \dots, b_g, c_1, \dots, c_{g-1}, e_0, \dots e_{n-1}$. 
	\end{theorem}
	
	We will use the following relations in the mapping class group in our proofs:
	\begin{proposition}
		\begin{enumerate}
			\item(Braid Relation)\cite[Proposition 3.11]{farbprimer} If $a$ and $b$ are simple closed curves intersecting once, then $T_aT_bT_a = T_bT_aT_b$
			\item(Lantern Relation)\cite[Proposition 5.1]{farbprimer} Let $S_0^4$ denote a sphere with 4 boundary components embedded in $S$. Denote the boundary curves by $b_1, \dots, b_4$ and let $x, y, z$ be non-isotopic separating simple closed curves in $S_0^4$ with $x \cap y = x \cap z = y \cap z = 2$. Then $T_x T_y T_z = T_{b_1} T_{b_2} T_{b_3} T_{b_4}$. 
		\end{enumerate}
	\end{proposition}
	
	We next record that Dehn twists about nonseparating curves normally generate $\Mod(g,n)$. 
	\begin{lemma}\label{normalclosure}
		Let $(g,n) \geq (1,0)$. Let $a$ be a nonseparating simple closed curve on a topological surface $\Sigma_{g,n}$, and let $T_a \in \Mod(g,n)$ denote the Dehn twist about $a$. Then the normal closure of $T_a$ is $\Mod(g,n)$. 
	\end{lemma}
	\begin{proof}
		For any simple closed curve $a$ and a mapping class $\phi \in \Mod(g,n)$, we have that $T_{\phi(a)} = \phi T_a \phi\inv$. Hence, all nonseparating simple closed curves (or more generally, simple closed curves of the same topological type) are conjugate. The result follows because the mapping class group is finitely generated by a set of nonseparating Dehn twists. 
	\end{proof}
	
	\subsection{Subgroups of the Mapping Class Group}
	In this section, we recall two distinguished subgroups of $\Mod(g,n)$ that will come up in our proofs. 
	
	\subsubsection{Torelli Subgroup}
	There is a natural representation $\Mod(g,n) \to \Sp(2g,\ZZ)$ given by the action of $\Mod(g,n)$ on the first homology of the surface without marked points, and it is known to be surjective. Its kernel is called the \textit{Torelli subgroup}, denoted $\shi_{g,n}$. Let $a' \subset S$ be a simple closed curve homologous to $a$.
	
	\begin{proposition}\label{torelligens}
		\begin{enumerate}
			\item\cite[Theorem 6.15 + $\varepsilon$]{farbprimer} When $g \geq 3$, $\shi_{g,n}$ is normally generated by $T_{a} T_{a'}^{-1}$.
			\item\cite[Corollary 1.4]{Putman} When $g = 2$, $\shi_{g,n}$ is normally generated by Dehn twists about separating curves wich partition $S$ into two genus one surfaces with one boundary component, one with $p$ marked points and one with $q$ marked points, with $p+q=n$, for all such partitions of $n$. 
		\end{enumerate}
	\end{proposition}
	
	\subsection{Point-Pushing Subgroup}
	There are natural "forget a point" maps $\Mod(g,n+1) \to \Mod(g,n).$ The kernel of this map is called the \textit{point-pushing subgroup}. There is a short exact sequence
	$$1 \to \pi_1(\Sigma_{g,n}) \to \Mod(g,n+1) \to \Mod(g,n) \to 1,$$ called the \textit{Birman exact sequence}. It is naturally identified as the homotopy exact sequence for the fibration $\shm_{g,n+1} \to \shm_{g,n}$, and by it we realize $\pi_1(\Sigma_{g,n})$ as the point pushing subgroup. 
	
	We abuse notation and also call the kernel of the composition $\Mod(g,n) \to \Mod(g,0)$ the point pushing subgroup, denoted $P_{g,n}$. 
	
	\begin{proposition}\cite[Fact 4.7]{farbprimer}\label{ptpushkernel}
		$P_{g,n}$ is normally generated by the "point pushing" maps $T_{e_i} T_{e_{i+1}}^{-1}$ for $i = 0, \dots, n-1$.
	\end{proposition}
	
	\subsubsection{Stabilizer Subgroups and Mapping Class Groups of Subsurfaces}
	If $R \subset S$ is a subsurface, we can extend a mapping class by the identity to obtain a homomorphism $\Mod(R) \to \Mod(S)$. Given an isotopy class of a simple closed curve $c$, we may consider the subgroup $\Stab(c)$ of mapping classes which preserve $c$ pointwise. Clearly, if $R$ is the complement of a regular neighborhood of $c$, $\Mod(R)$ maps to $\Stab(c)$, and the kernel identifies the two boundary twists with $T_c$. 
	
	On the other hand, keeping $R$ to be the complement of a neighborhood of $c$, so that $R$ has genus $g-1$ and $2$ boundary components, we may cap the boundary components with punctured disks, obtaining the short exact sequence
	$$1 \to \{T_{\partial_1}, T_{\partial_2}\} \to \Mod(R) \to \Mod(g-1,n+2)_{p_1, p_2} \to 1.$$
	
	Here the subscript denotes that we demand the mapping classes preserve the two new punctures individually. This will not make a difference for most of our arguments. When workign with moduli, we always take the fundamental group of $\shm_{g,n}$ to be the pure mapping class group, which already preserves all marked points individually. 
	
	Deleting the curve $c$ gives a map from $\Stab(c)$ to $\Mod(g-1,n+2)_{p_1, p_2}$. If we let $\Stab(c)^+$ denote the subgroup of $\Mod(g,n)$ stabilizing $c$ and preserving its orientation, the cutting map $\zeta: \Stab(c)^+ \to \Mod(g-1,n+2)_{p_1, p_2}$ is surjective. In particular, if $Cap:\Mod(R) \to \Mod(g-1,n+2)_{p_1, p_2}$ is the capping homomorphism, $p:\Mod(g-1,n+2) \to \Mod(g-1,n+2)_{p_1, p_2}$ the quotient by the involution of the two points, and $\iota: \Mod(R) \to \Stab(c)^+$ the inclusion, we have $\zeta \circ \iota = p \circ Cap$. 
	
	We have the following short exact sequence:
	$$1 \to \langle T_c \to \rangle \Stab(c)^+ \to \Mod(g-1,n+2)_{p_1, p_2} \to 1.$$ 
	
	Later, we will compare this to the short exact sequence for a punctured tubular neighborhood of the boundary divisor $\Delta_0$ in $\shm_{g,n}$. For now, we work purely group-theoretically. Suppose $\rho$ is a representation of $\Mod(g,n)$ such that $\rho|_{\Stab(c)^+}$ preserves a filtration $W_\bullet$. Suppose further that there is an induced representation $\bar{\rho}: \Mod(g-1,n+2)_{p_1, p_2} \to GL(\gr_W V)$. Let $\gr \rho: \GL(V) \to \GL(\gr_W V)$ denote the associated graded. If $d \subset S$ is disjoint from $c$, we may consider $d$ as a curve in $R$, hence $\gr \rho (T_d) = \bar{\rho} \circ \zeta (T_d)$. After restricting to $\Mod(R)$, we have $\gr \rho \circ \iota = \bar{\rho} \circ Cap$. From this, we can gather information about the image of a nonseparating Dehn twist, particularly the eigenvalues, by restricting to the subsurface $R$.

	\subsection{Abelianization of the Mapping Class Group}
	The most important input to the classification of low-dimensional representations of the mapping class group is the fact, originally due to Mumford and expanded by others, that for $g \geq 2$, $\Mod(g,n)$ has finite abelianization. 
	\begin{theorem}\cite[Theorem 2.3]{korkmaz}\label{finiteab}
		The abelianization of the mapping class group $\Mod(g,n)$ is a quotient of $\ZZ/10\ZZ$ for $g \geq 2$, $n \geq 0$. For $(g,n) \geq (3,0)$, it is trivial. 
	\end{theorem}
	
	A related fact, which we will use repeately, is about the normal generation of the commutator subgroup:
	
	\begin{proposition}\cite[Theorem 2.4]{korkmaz}
		When $ g\geq 2$, the commutator subgroup of $\Mod(g,n)$ is normall generated by $T_aT_b^{-1}$, where $a$ and $b$ are nonseparating simple closed curves which intersect at a point. 
	\end{proposition}
	
	In particular, for any representation $\rho$ of $\Mod(g,n)$ with $a$ and $b$ as above, if $\rho(T_a) = \rho(T_b)$, $\rho$ annihilates the commutator subgroup of $\Mod(g,n)$.

\section{Classification of Representations of the Mapping Class Group}
	\subsection{Preliminary Results}
		We will use the following criterion of Korkmaz to determine when a representation of the mapping class group is trivial. 

		We will use the following criteria of Korkmaz to find invariant subspaces for a representation of $\Mod(g,n)$.
	
		\begin{lemma}\cite[Lemma 4.3]{korkmaz}
			Let $a,b$ be two nonseparating simple closed curves on $S$ intersecting once, and suppose $\lambda$ is an eigenvalue of both $\rho(T_a)$ and $\rho(T_b)$. If the $\lambda$-eigenspaces of $\rho(T_a)$ and $\rho(T_b)$ are equal, then that eigenspace is invariant under $\Mod(g,n)$. 
		\end{lemma}

		In fact, Korkmaz's proof of \Cref{preserveespaces} shows a little more.
		
		\begin{lemma}\label{preserveespaces}
			Let $\rho: \Mod(g,n) \to \GL(r,\CC)$ be a representation with $g \geq 2$ and let $a$ and $b$ be two nonseparating simple closed curves which intersect at a point. Suppose $V_a$ is a subspace preserved by $\rho(T_a)$ and by the image of every mapping class commuting with $\rho(T_a)$. If $\rho(T_b)$ preserves $V_a$, then $V_a$ is $\Mod(g,n)$-invariant. 
		\end{lemma} 
			\begin{proof}
				Let $\phi \in \Mod(g,n)$ be a mapping class sending the curve $a$ to a curve $x$. Define $V_x := \rho(\phi)V_a$; this is well-defined since mapping classes commuting with $T_a$ preserve the isotopy class of $a$. Since $T_aT_bTa$ sends $a$ to $b$, while $\rho(T_bT_aT_b)V_a = V_b$ we have $V_b = V_a$. If $x$ and $y$ other are nonseparating simple closed curves intersecting at a point, there is a mapping class $\phi$ sending the pair $(a,b)$ to $(x,y)$, hence $V_x = V_y$. Since any two nonseparating simple closed curves may be connected by a chain of nonseparating simple closed curves, and there exists a diffeomorphism mapping the pair $(a,b)$ to any link of this chain, it follows that $V_c = V_a$ for every nonseparating simple closed curve $c$. 
			\end{proof}
			
		A common application of this lemma will be that if the image of a nonseparating dehn twist has multiple eigenvalues, then the braid relation often forces one of those eigenspaces to be $\Mod(g,n)$-invariant. More generally, if a subspace $W$ is preserved by the image of the centralizer of $\rho(T_a)$ and by $\rho(T_b)$, then $W$ is $\Mod(g,n)$-invariant. 
		
	\subsection{Triviality of Low Dimensional Representations}

		Our proof begins with the base case of $g = 2$ and $r \leq 3$. For this base case the proof is essentially the same as Korkmaz, though we may eliminate some casework using Hodge theory. Here and throughout, we use the notation $E_\lambda^c$ to mean the $\lambda$ eigenspace of $\rho(T_c)$, i.e. $\ker(\rho(T_c) - \lambda I).$ We will occastionally use $(E_\lambda^c)_{gen}$ to refer to the generalized eigenspace.

		\begin{proposition}\label{g2finiteimage}
			Let $\rho: Mod(2,n) \to \GL(r, \CC)$ be a representation with $r \leq 3$. Then the image of $\rho$ is a quotient of $\ZZ/10\ZZ$. 
		\end{proposition}
		\begin{proof}
			We may deform $\rho$ so that it underlies a $\CC$-variation of Hodge structure and denote by $k$ the weight of the Hodge structure. To $\rho$ we may associate a period map on $\shm_{g,n}$, and the local monodromy of this period map around the divisor $\Delta_0$ is given by $\rho(T_c)$, where $c$ is a nonseparating simple closed curve. By \ref{qudensity}, we may assume $\rho(T_c)$ is quasi-unipotent.If $r = 1$, we are done by \Cref{finiteab}. Write $\rho(T_c) = SU$ where $S$ is has finite order and $U = I+N$, with $N = 0$ having nilpotency order $k+1$. When $r = 2$, only two orders of nilpotency are possible: 0 or 1. Consider two nonseparating simple closed curves  $a$ and $b$, both disjoint from $c$, which intersect once. We can arrange so that $\rho(T_a)$ and $\rho(T_b)$ are both upper-triangular with respect to the Jordan basis of $\rho(T_c)$. 
			
			We will consider each possible rank of $N$, in each case proving the statement for $r=2$ and $r=3$ individually. 
			\textbf{$N = 0$. } 
			
			\textit{(a) $r=2$.} There are two subcases. If $\rho(T_c)$ has a unique eigenvalue, then since all Dehn twists about nonseparating curves are conjugate, $\rho(T_a)$ and $\rho(T_b)$ must be of the same form:
			\[\rho(T_a) \sim \rho(T_b) \sim \begin{bmatrix}
				\lambda & *\\ & \lambda
			\end{bmatrix}.\] The braid relation then implies $\rho(T_a) = \rho(T_b)$, so that $\rho$ annihilates the commutator subgroup of $\Mod(2,n)$.
			
			 The other possibility is that $\rho(T_c)$ has two distinct eigenvalues, but in this case $\rho(T_a)$ and $\rho(T_b)$ must both be semisimple, hence they commute, and again the braid relation forces their images to be equal. For $\rho(T_c)$ to have a nontrival unipotent part, it has a unique eigenvalue. Consequentially, $\rho(T_a)$ and $\rho(T_b)$, which are both conjugate to $\rho(T_c)$, must commute with each other. By the braid relation, $\rho(T_a) = \rho(T_b)$, so $\rho$ annihilates the commutator subgroup of $\Mod(2,n)$. 

			\textit{(b) $r=3$.} Again, consider the nilpotent part of $\rho(T_c)$ and its order of nilpotency: 0, 1, or 2. If $0$, the image of $\rho$ is abelian by similar arguments as in the case $r=2$. Indeed, if there are three unique eigenvalues, we are done by the same argument verbatim. Similarly if there is a unique eigenvalue, consider a nonseparating simple closed curve $d$ intersecting $c$ once. Fixing a basis where $\rho(T_c) = \lambda I$, $\rho(T_c)$ centralizes $\rho(T_d)$, and the braid relation implies $\rho(T_c) = \rho(T_d)$. Repeating this throughout a chain of nonseparating simple closed curves shows that every generator of $\Mod(2,n)$ maps to $\lambda I$.
			
			If there are two unique eigenvalues $\lambda$ and $\mu$ with $\lambda$ having multiplicity two, let $\chi: \Mod(2,n) \to \ZZ/10\ZZ$ be the character corresponding to $\lambda$ and twist $\rho$ by $\chi\inv$, so that we may assume the eigenvalues are $1$ and $\mu$. Again consider a nonseparating simple closed curve $d$ intersecting $c$ once. If the $\lambda-$eigenspaces for $\rho(T_c)$ and $\rho(T_d)$ are equal, then we are done by \Cref{preserveespaces}. Otherwise, they overlap in one dimension. Pick a basis for $\CC^3$ with $v_1 \in E_\lambda^{c} \setminus E_\lambda^d, v_2 \in E_\lambda^c \cap E_\lambda^d,$ and $v_3 \in E_\mu^c$. 

			Now consider two nonseparating simple closed curves $a$ and $b$, both disjoint from $c$ and $d$, which intersect at a point. By the braid relation, the diagonal entries of $\rho(T_a)$ and $\rho(T_b)$ are equal, so there are three possibilities: $\rho(T_a) = \rho(T_b) = \rho (T_c)$, in which case the image is abelian and we are finished, or either: 
			$$\rho(T_a) = \begin{bmatrix}
				\lambda & a & \\
				 & \mu & \\
				 & & \lambda \\
			\end{bmatrix}, \qquad \text{or} \qquad \rho(T_a) = \begin{bmatrix}
				\mu & a & \\
				& \lambda & \\
				& & \lambda \\
			\end{bmatrix}.$$ $\rho(T_b)$ has the same form, but with a different off-diagonal term $b$. In both cases the braid relation implies $(\mu^2 - \lambda\mu + \lambda^2)(a-b)=0$. Now let $e$ be a nonseparating simple closed curve disjoint from $a$ and $c$, intersecting $b$ and $d$ once each. We may write $\rho(T_e)$ in the form above, and since $[\rho(T_e),\rho(T_a)] = 0$, we have $(\lambda-\mu)(a-e)$ = 0, where $ a$ and $e$ are their unique off-diagonal entries. Similarly the braid relation with $\rho(T_b)$ gives $(\mu^2 - \lambda\mu + \lambda^2)(b-e) = 0$. If $\mu^2 - \lambda\mu + \lambda^2 = 0$, we must have $a=e$, so that $\rho(T_e) = \rho(T_a)$. But now the braid relation with $\rho(T_d)$, combined with the fact that $\rho(T_d)$ commutes with $\rho(T_a) = \rho(T_e)$, means $\rho(T_d) = \rho(T_a)$ as well. We arrive at a number of contradictions, but in particular we find that the image of $\rho$ is abelian.
			
			\textbf{$\rk N=1$.} If the order of nilpotency is one, $\rho(T_c)$ has one or two distinct eigenvalues. When there are two eigenvalues, the eigenspace for the eigenvalue of lesser multiplicity is $\Mod(2,n)-$invariant, hence $\rho$ is reducible and we are done by the rank 1 and 2 cases.
			
			In the two cases left to consider, $\rho(T_c)$ has a unique eigenvalue $\lambda$. Taking $a_1, b_2, a_2, b_2, c$ to be the Humphries generators of $\Mod(2)$ with $[c] \equiv [a_1] - [a_2]$, the chain relation $(T_{a_1} T_{b_1} T_c)^4 = T_{a_2}^2$ shows that $\lambda^{12} =\lambda^2$. Without loss of generality, we may replace $\rho$ with the representation $\rho \otimes \lambda^{-1}$, for which the images of nonseparating Dehn twists are unipotent. 

			Write $\rho(T_c) = I + N$ with N nilpotent. The order of nilpotency is either 2 or 3, corresponding to the rank of $N$ being one or two, respectively. If $\rk N = 1$, a direct calculation with the braid relation shows that $\rho(T_a) = \rho(T_b)$. 
			
			\textbf{$\rk N=2$.} If $\rk N = 2$, then the centralizer of $\rho(T_c)$ consists of polynomials in $N$,and this centralizer is abelain. Since $\rho(T_a)$ and $\rho(T_b)$ commute with $\rho(T_c)$, they lie in this centralizer, and the braid relation forces them to be equal.
		\end{proof}

		We may now prove \Cref{nolowrankreps}. We will prove this by induction, deforming $\rho$ to a local system which underlies a polarizable complex variation of Hodge structures. The result will quickly follow by induction when the local monodromy around the divisor $\Delta_0$ has a nontrivial unipotent part. In the inductive argument, we will make use of the following criterion of Korkmaz. Korkmaz proves uniformly in $g$ as a consequence of \Cref{2g-1trivial}. In our context of proving \Cref{2g-1trivial}, for a fixed $g$, it will hold for all genera up to $g-1$. 
		
		\begin{lemma}\Cite[Lemma 7.1]{korkmaz}\label{flag}
			Let $\rho: \Mod(g,n) \to \GL(r,\CC)$ be a representation, and assume \Cref{2g-1trivial} holds for all $g \in \{2, \dots, g\}$. Suppose there is a $\Mod(g,n)$ invariant filtration $W_\bullet$ whose graded pieces all have dimension at most $2g-1$. Then $\rho$ annihilates the commutator subgroup of $\Mod(g,n)$.
		\end{lemma}
			\begin{proof}
				The filtration $W_\bullet$ on $V = \CC^r$ gives a representation $\gr_W \rho: \Mod(g,n) \to \GL(\gr_W V)$. If each graded piece has dimension at most $2g-1$, then \Cref{2g-1trivial} applies, so that each of these representations annihilate the commutator subgroup of $\Mod(g,n)$. It follows that $\rho([\Mod(g,n), \Mod(g,n)])$ is contained in a solvable group, which, since the commutator subgroup is perfect, forces $\rho$ itself to annihilate the commutator subgroup of $\Mod(g,n)$.
			\end{proof}

		In order for $\rho$ to underlie a polarizable complex variation of Hodge structures with a nontrivial monodromy weight filtration, we must dispence with the case that nonseparating Dehn twists are semisimple.
		
		\begin{lemma}\label{nosstwists}
			Let $\rho: \Mod(g,n) \to \GL(r,\CC)$ be a representation with $(g,n) \geq (3,0)$ in lexicographic order, and $r \leq 2g-1$, such that the image of a Dehn twist about a nonseparating simple closed curve $c$ is semisimple. Then the image of $\rho$ is trivial. 
		\end{lemma}
			\begin{proof}
				We first consider when $\rho(T_c)$ has at least three eigenvalues. Then one may break up the filtration $0 \subset E_{\lambda_1}^c \subset E_{\lambda_1}^c \oplus E_{\lambda_2}^c \subset \dots \subset \bigoplus E_{\lambda_i}^c$ of generalized eigenspaces into two pieces whose graded dimensions are both at least two. Call that filtration $W_\bullet$. For example if $E_{\lambda_1}^c$ and $E_{\lambda_2}^c$ are both one-dimensional one can take $W_1 = E_{\lambda_1}^c \oplus E_{\lambda_2}^c$ and $W_2 = \CC^r$. 
				
				Let $R$ be the complement of a neighborhood of $c$ and consider $\rho|_{\Mod(R)}$. Every element of $\Mod(R)$ preserves $W_\bullet$. The representations on the graded pieces are representations of $\Mod(R)$ that satisfy the inductive hypothesis, so they are trivial when $g \geq 4$ and have cyclic image when $g = 3$. Since the commutator subgroup of $\Mod(R)$ is perfect, the image of $\Mod(R)$ is trivial by \Cref{flag}. Since the normal closure of the commutatuor of $\Mod(R)$ in $\Mod(g,n)$ is $\Mod(g,n)$, the image of the commutator subgroup of $\Mod(g,n)$, is trivial. Since $g \geq 3$, this implies $\rho$ is trivial.
				
				Now we return to the situation where $\rho(T_c)$ has only one or two eigenvalues. If one of the (generalized) eigenspaces of $\rho(T_c)$ is such that both $E_\lambda^c$ and $\CC^r/E_\lambda^c$ have dimension at most $r-2$, the same inductive argument as above applies. So we are only left with the cases where $\rho(T_c) = \lambda I$ or $\rho(T_c) = \text{diag}(\lambda, \dots, \lambda, \mu)$. In the former case, of course, conjugacy of all nonseparating Dehn twists implies that every nonseparating Dehn twist maps to $\lambda I$. 
				
				We are left to consider when $\rho(T_c)$ has two distinct eigenvalues, in which case the filtration $W_\bullet$ above coming from the Jordan filtration does not allow for an induction. or $\rho(T_c) = \text{diag}(\lambda, \dots, \lambda, \mu).$
				Consider a nonseparating simple closed curve $d$ intersecting $c$ once. If the $\lambda-$eigenspaces of $c$ and $d$ are equal, write $\rho(T_c) = \lambda I+c$ and $\rho(T_d) = \lambda I+d$, where $\ker c = \ker d = E_\lambda^c$. Since $\img d \subset \ker d = \ker c$, $dc = 0$, and likewise $dc=0$. On the other hand the braid relation expands to $\lambda I+2c+d = \lambda I+c+2d$, hence$c=d$. Therefore $\rho$ annihilates the commutator subgroup of $\Mod(g,n)$, so the image of $\rho$ is trivial. 
				
				We are left with the case where the $\lambda$-eigenspaces of $\rho(T_c)$ and $\rho(T_d)$ are different. In this case, the filtration $0 \subset E_\lambda^c \cap E_\lambda^d \subset \CC^{r}$ together with the inductive argument above shows that the image of $\rho$ is trivial.  
			\end{proof}
		
		\begin{theorem}\label{2g-1trivial}
			Let $\rho: \Mod(g,n) \to \GL(r, \CC)$ be a representation with $(g, n) \geq (2, 0)$ and $r \leq 2g-1$. If $g=2$ then the image of $\rho$ is a quotient of $\ZZ/10\ZZ$, while if $g \geq 3$, $\rho$ is trivial. 
		\end{theorem}
			\begin{proof}
				The case when $g = 2$ was done in \Cref{g2finiteimage}. We now proceed by induction on $g$. We assume $\rho$ underlies a $\CC$-variation of Hodge structures with quasi-unipotent local monoromy around the divisor $\Delta_0$, hence gives a period map $f: \shm_{g,n} \to D/\Gamma$ to some period domain. Let $a$ be a nonseparating simple closed curve on the underlying genus $g$ surface $S$. Then the monodromy around the boundary component $\Delta_0 \subset \bar{\shm}_{g,n}$ is given by $\rho(T_a)$, which we assume to be quasi-unipotent, with unipotent part equal to $I+N$. We may assume $N$ is nontrivial, since otherwise we are done by \Cref{nosstwists}. Then $N$ induces the monodromy weight filtration $W_\bullet$ on $\CC^r$. 
				
				Let $R \subset S$ be the complement of a regular neighborhood of $a$. Then the genus of $R$ is $g-1$, and there is a natural inclusion $\Mod(R) \to \Mod(g,n)$ given by extension by the identity. In particular, the image of $\Mod(R)$ commutes with $T_a$, so $\rho(\Mod(R))$ preserves $W_\bullet$. Since $W_\bullet$ is nontrivial, the dimension of each of these graded pieces is at most $r-2 \leq 2g-3 = 2(g-1) -1$. By \Cref{flag}, $\rho$ annihilates the commutator subgroup of $\Mod(R)$. The normal closure of the commutator subgroup of $\Mod(R)$ in $\Mod(g,n)$ is $\Mod(g,n)$, hence $\rho$ factors through $\Mod(g,n)^{ab} = \{1\}$.  
			\end{proof}
			
		For maps $\shm_{g,n} \to \sha_h$ with $h < g$, we can give a geometric interpretation of this argument: $\shm_{2,n}$ cannot map to $\sha_1$ becuase $H^1(\shm_{g,n}, \CC) = 0$. Now if $\shm_{g,n}$ maps to $\shm_g$ with nontrial unipotent monodromy around $\Delta_0$, there is an induced map $\shm_{g-1,n+2} \to \sha_{h-r}$ for some $r >0$. Inductively, this map must be constant, and the lift of any Dehn twist from $\Mod(g-1,n+2)$ to $\Mod(g,n)$ will then lie in a solvable subgroup of $\Mod(g,n)$, showing that the original representation was trivial. This argument will be incorporated into the inductive argument of \Cref{mga2g-1} below, where we justify that the local monodromy around $\Delta_0$ is unipotent and handle the case where $\Delta_0$ is contracted to $\sha_0$.
	\subsection{Rank $2g$ Representations of the Mapping Class Group}
		We now begin the classification of rank $2g$ representations of the mapping class group. We proceed by induction on $g$ and $n$. We establish the case $g=2$ explicitly, and our argument is in the same spirit as Korkmaz's. For $g \geq 3$, we may use \Cref{MocDeform} to assume $\rho$ underlies a complex variation of Hodge structure. We then consider the local monodromy around the divisor $\Delta_0 \subset \bar{\shm}_{g,n}$. An analysis of a punctured tubular neighborhood around $\Delta_0$ lets us determine the image of $\rho$ on most of $\Mod(g,n)$, and a calculation with the braid relation lets us determine the rest of the image. We conclude that $\rho$ is the standard symplectic representation, which is rigid, so our initial deformation was trivial. 
		\subsubsection{The Base Case: $g=2$}
		The following argument occurs often enough in the rest of the paper that we write it in an explicit basis. 
		\begin{lemma}\label{usefulcalcs}
			\begin{enumerate}
				\item 
				Let $$U=
				\begin{bmatrix}
					\lambda & 1 \\
					0 & \lambda
				\end{bmatrix}.$$ A matrix $V$ commuting with $U$ must be of the form $$V=
				\begin{bmatrix}
				a & b \\
				0 & a
				\end{bmatrix}.
				$$
				\item The braid relation in rank 2: If $$W = \begin{bmatrix}
					\lambda & 0\\ -\lambda^2 &\lambda\\
				\end{bmatrix},$$ then $UWU=WUW$. If $WVW=VWV$, then $U=V$.
			\end{enumerate}
		\end{lemma}
			\begin{proof}
				Item (1) and the first part of (2) are direct computations:
				\[
				\begin{bmatrix}
					a & b\\
					c & d
				\end{bmatrix}
				\begin{bmatrix}
					\lambda & 1\\
					0 & \lambda
				\end{bmatrix}
				=
				\begin{bmatrix}
					\lambda & 1\\
					0 & \lambda
				\end{bmatrix}
				\begin{bmatrix}
					a & b\\
					c & d
				\end{bmatrix}.
				\]
				Thus
				\[
				\begin{bmatrix}
					a\lambda & a+b\lambda\\
					c\lambda & c+d\lambda
				\end{bmatrix}
				=
				\begin{bmatrix}
					a\lambda+c & b\lambda+d\\
					c\lambda & d\lambda
				\end{bmatrix}.
				\]
				Hence $c = 0$ and $a = d$.
				Similarly \[ UVU
				=
				\begin{bmatrix}
					0 & \lambda\\
					-\lambda^3 & \lambda^2
				\end{bmatrix}
				\begin{bmatrix}
					\lambda & 1\\
					0 & \lambda
				\end{bmatrix}
				=
				\begin{bmatrix}
					0 & \lambda^2\\
					-\lambda^4 & 0
				\end{bmatrix},
				 \quad \text{and} \quad
				VUV
				=
				\begin{bmatrix}
					\lambda^2 & \lambda\\
					-\lambda^3 & 0
				\end{bmatrix}
				\begin{bmatrix}
					\lambda & 0\\
					-\lambda^2 & \lambda
				\end{bmatrix}
				=
				\begin{bmatrix}
					0 & \lambda^2\\
					-\lambda^4 & 0
				\end{bmatrix}.\]
				
				 For the second part of 2, note that if $W$ and $V$ satisfy the braid relation, they are conugate, hence $a = \lambda$. From there, the computation $VWV=WVW$ shows $b = 1$. 
			\end{proof}

		Before proving the uniqueness of the symplectic representation in genus 2 in full, we separately treat one of the more involved cases of the local monodromy around $\Delta_0$. 
		
		\begin{lemma}\label{no22jordanform}
			Let $\rho: \Mod(2,n) \to \GL(4,\CC)$ be a representation which underlies a polarizable complex variation of Hodge structures of weight one. If the rank of the nilpotent part of the monodromy around $\Delta_0$ is two, then the image is a quotient of $\ZZ/10\ZZ$. 
		\end{lemma}
			\begin{proof}
				Let $a$ and $b$ be nonseparating simple closed curves intersecting at a point. We are ruling out a local monodromy matrix of the form $$\rho(T_a) = \begin{bmatrix}
					\lambda &1 & & \\
					& \lambda & & & \\
					& & \mu & 1\\
					& & & \mu
				\end{bmatrix}.$$ 
				First suppose $\lambda \neq \mu$. Then the $\lambda$ eigenspaces of $\rho(T_a)$ and $\rho(T_b)$ are one dimensional. If there is any nontrivial overlap between the two, then there is an invariant line under $\rho$, and we are done by the classification of lower dimensional representations of $\Mod(2,n)$. If the $\lambda$ and $\mu$ eigenspaces are disjoint, we first reduce to the case $n=0$, i.e. we get rid of points. If $a'$ is a nonseparating simple closed curve isotopic to $a$ such that $a$ and $a'$ bound a point (in particular, the two curves are disjoint), we may write $\rho(T_{a'})$ to be block diagonal with respect to the Jordan basis for $\rho(T_a)$ since it preserves the generalized eigenspaces of $\rho(T_a)$. Then \Cref{usefulcalcs} shows that $\rho(T_a) = \rho(T_{a'})$ so that $\rho(T_a T_{a'}\inv) = I$. 
				
				Let $a,b,c,d,e$ be a chain of nonseparating simple closed curve whoses twists generate $\Mod(2)$. The chain relation $(T_a T_b T_c)^4 = T_e^2$ implies that $\det(\rho(T_a)) = (\lambda \mu)^2$ is a 10th root of unity. The classification of 4 dimensional representations of the braid group $B_3$ given in \Cite[Proposition 2.6]{b3reps} gives explicit bases  $\rho(T_a)$ and $\rho(T_b)$, under which $\rho(T_a T_b)^6 = (\lambda \mu)^6$. This element represents the image of a separating twist, hence the image of the Torelli subgroup has order at most 10. Moreover, since $\rho(\mathcal{I}_2)$ is diagonal, it commutes with $\rho(\Mod(2))$, hence its $\lambda$ and $\mu$ eigenspaces are invariant, so by \Cref{2g-1trivial} $\rho$ itself factors through $\Mod(2)_{ab}$.

				Now assume $\lambda = \mu$. By multiplying by the matrix $\lambda\inv I$, we may assume $\rho(T_a)$ and $\rho(T_b)$ are unipotent. The $1-$eigenspaces of each are two dimensional, so they have 0,1, or 2 dimensional overlap. If they overlap in two dimensions, the eigenspaces are equal, hence invariant by \Cref{preserveespaces}, and we are done. Now assume the two eigenspaces are disjoint. In block form we may write
				$$\rho(T_a) = \begin{bmatrix}
					I_2 & I_2\\ 0 & I_2\\
				\end{bmatrix}, \qquad \rho(T_b) = \begin{bmatrix}
				I_2 & 0 \\ B & I_2\\
				\end{bmatrix},$$ with $B \in \GL(2,\CC)$. The braid relation gives $I+2B = B$, so $B=-I$. Then $(\rho(T_a)\rho(T_b))^6 = I$, so $\rho(\mathcal{I}_2) = I$, and $\rho$ factors through the surjection to $\Sp(4,\ZZ)$. Under this surjection, $T_a$ maps to a unipotent matrix with rank one. This is a contradiction because the only nontrivial representation $\Sp(4) \to \GL(4)$ is the standard representation. 
				
				We are left with the case where the $\lambda-$eigenspaces of $\rho(T_a)$ and $\rho(T_b)$ overlap in one dimension. Let $v \in E_\lambda^a \cap E_\lambda^b$. Because $\rho(T_d)$ and $\rho(T_e)$ preserve the $\lambda-$eigenspaces of both $\rho(T_a)$ and $\rho(T_b)$, they preserve the line spanned by $v$. To be precise, they both send $v$ to $(\lambda v)$. Hence $v \in E_\lambda^a \cap E_\lambda^e$. We claim this is an equality. If not, the two eigenspaces are equal. Since $\rho(T_b)$ preserves $E_\lambda^e = E_\lambda^a$, the braid relation implies $E_\lambda^a = E_\lambda^b$: If $N_a$ and $N_b$ are the nilpotent parts of $\rho(T_a)$ and $\rho(T_b)$, the braid relation implies $N_a +N_aN_bN_a = N_b + N_bN_aN_b$. Letting $w \in E_\lambda^a$, $N_a w = 0$, so $(N_b+N_bN_aN_b)w = 0$. Since $\rho(T_b)$ preserves $E_\lambda^a$, $N_a(N_b w) = 0$, so $N_b w = 0$, hence $w \in E_\lambda^b$. It follows that the line spanned by $v$ is invariant under $\rho$ --- in fact, $\rho$ sends $v$ to $\lambda v$ with $\lambda$ a $10$th root of unity, so by the theorem of the fixed part, $\rho$ has an isotrivial isogeny factor, hence the whole representation has image factoring through $\ZZ/10\ZZ$. 
			\end{proof}
			
			As a consequence of this, we are reduced to classifying representations of $\Mod(2,n)$ for which the image of a nonseparating Dehn twist has a unique eigenvalue $\lambda$, which, by the chain relation, is a 10th root of unity. 
			
		\begin{corollary}\label{only2eval}
			Let $\rho: \Mod(2,n) \to \GL(4,\CC)$ be a representation which underlies a polarizable complex variation of Hodge structure on $\shm_{2,n}$ with quasi-unipotent monodromy around $\Delta_0$. If $\rho$ has infinite image, then the image of a nonseparating Dehn twist has a unique eigenvalue, a 10th root of unity, and the $\lambda$-eigenspace is 3-dimensional. 
		\end{corollary}
			\begin{proof}
				If there are three or four distinct eigenvalues, one of the eigenspaces is one dimensional, hence $\Mod(2,n)$-invariant. Indeed, for example when $\rho(T_a)$ has three distinct eigenvalues $\lambda, \mu,$ and $\tau$, consider $c$ and $d$ two nonseparating simple closed curves disjoint from $a$ which intersect at a point. Using that $\rho(T_c)$ and $\rho(T_d)$ are conjugate to $\rho(T_a)$ and preserve its eigenspaces, we have

				\[\rho(T_c) \sim \rho(T_d) \sim \begin{bmatrix}
					\lambda & * & & \\
					*&\lambda & & \\
					& & \mu & \\
					 & & & \tau
				\end{bmatrix}.\] The braid relation implies the $\mu$ and $\tau$ eigenspaces of $\rho(T_c)$ and $\rho(T_d)$ are equal, and \Cref{preserveespaces} implies they are $\Mod(2,n)$ invariant. It follows that $\rho$ is reducible, hence it has finite image by our classification of lower rank representations.
				
				If there are two distinct eigenvalues, then we are either in a case handled by the previous lemma, or one of the eigenspaces is one dimensional, or in the case where the Jordan form of $\rho(T_a)$ is $(2,1,1)$ which we now treat. 
				
				Let $d$ and $e$ be two nonseparating simple closed curves intersecting at a point, both disjoint from $a$. Since $\rho(T_d)$ and $\rho(T_e)$ are conjugate to $\rho(T_a)$ and both commute with $\rho(T_a)$, in the Jordan basis for $\rho(T_a)$ we may write $\rho(T_d)$ and $\rho(T_e)$ to be upper triangular. Without loss of generality, suppose that the two eigenvalues are $1$ and $\lambda$, with $\lambda$ having lesser multiplicity. If in that basis, $E_\lambda^d = E_\lambda^e$, we are done, since $\rho$ then has an invariant line. If not, then we may write
				$$\rho(T_d) = \begin{bmatrix}
					1 & b & c & \\ & 1 & & \\ & d & \lambda & \\ & & & 1
				\end{bmatrix}, \qquad \rho(T_e) = \begin{bmatrix}
				1 & x & y &\\ & 1 & &\\ & z & \lambda & \\ & & & 1\\
				\end{bmatrix}.$$
				
				The braid relation gives 
				$$\rho(T_b T_c T_b - T_c T_b T_c) = \begin{bmatrix}
					0 &
					b-x+\lambda(cd-y\delta)&
					(c-y)(\lambda^2-\lambda+1)&0
					\\
					0 &0 &0 &0 \\
					 0&
					(d-z)(\lambda^2-\lambda+1)&0
					 &0 \\
				0 &0 &0 & 0
				\end{bmatrix}.$$
				 Since $\lambda = \det(\rho(T_a)) \in \ZZ/10\ZZ$, the term $\lambda^2 - \lambda +1$ is nonzero, so $c=y$ and $d=z$, hence $b=x$ as well, so $\rho(T_d) = \rho(T_e)$ and $\rho$ annihilates the commutator subgroup of $\Mod(2,n)$.
			\end{proof}
		\begin{proposition}\label{sp4rigid}
			Let $\rho: Mod(2,n) \to \GL(4,\CC)$ be a representation. If the image of $\rho$ is infinite, then $\rho$ is conjugate to either $\std$ or $\std \otimes \chi$, where $\std$ is the standard representation of $\Sp(4,\ZZ)$ and $\chi: \Mod(2,n) \to \{\ZZ/10\ZZ \}$ is the a character. Otherwise, the image of $\rho$ is a quotient of $\ZZ/10\ZZ$. 
		\end{proposition}
			\begin{proof}
				We assume $\rho$ is irreducible, for if $\rho$ is reducible, \Cref{g2finiteimage} shows that the image of $\rho$ is a quotient of $\ZZ/10\ZZ$. Let $a, b, c, d,e \subset \Sigma_{2,n}$ be a set of nonseparating simple closed curves on $S$ in standard position. Let $a= a_1, a_2 \dots, a_n$ be simple closed curve homologous to $a$ such that $a_i$ and $a_{i+1}$ bound the $i$th marked point. Let $A,B,C,D,E, \{A_i\}$ denote their images under $\rho$. Since $d$ and $e$ are disjoint from $a$ and $b$, $D$ and $E$ preserve the generalized eigenspaces of $A$ and $B$.				
				
				By \Cref{only2eval}, we may assume $A$ has a unique eigenvalue $\lambda$, which is a 10th root of unity, and that the $\lambda$-eigenspace has codimension one. By twisting $\rho$ by the character $\lambda\inv$, we may assume $A$ is uniotent. Write $A = I+a$ and	$B = I +b$, with $a$ and $b$ having rank one, the braid relation $ABA=BAB$ gives $a+aba = b+bab.$ If $a$ and $b$ commute, we have $a=b$ and $\rho$ annihilates the commutator subgroup. Otherwise, we have $aba \in \img(a)$ and $bab \in \img b$, hence $aba = \eta a$ and $bab= \mu b$ for $\eta, \mu \neq 0$. By associativity, $a \mu b = a(bab) = (aba)b = \eta a b$, we have $\mu = \eta$. From the braid relation, we get that $\eta = -1$. As a consequence we see that $a+aba = 0$.
				
				Fix $i \in 1, \dots, n$ and write $A_i = I + a'$ with $\rk a = \rk a' = 1$ and $[a,a'] = 0$. Two rank one nilpotent operators that commute have images on the same line or disjoint from each other. In the former case, we have $a' = \eta a$, and the analysis of the previous paragraph says $a = a+aba = b = a'+a'ba' = a'$, so that $A = A_i$, hence $\rho$ annihilates the point pushing subgroup. The latter case turns out to be impossible, since $A$ and $A'$ preserve $E_1^d \cap E_1^e$, hence act on a 2 dimensional space, where two commuting nilpotent operators which satisfy the braid relation with $b$ must be equal. 
				
				The work in the previous paragraph lets us reduce to a representation $\rho: \Mod(2) \to \GL(4,\CC)$. A regular neighborhood of $a \cup b$ is a torus with one boundary curve $\delta$. From the braid relation, $T_\delta = (T_a T_b)^6$, hence $\rho(T_\delta) = (AB)^6 = I+12(a+aba) = I$. Thus $\rho \otimes \chi$ factors through the surjection $\tau: \Mod(2) \to \Sp(4,\ZZ)$. So $\rho \otimes \chi = \varphi \circ \tau \circ \pi$, with $\varphi: \Sp(4,\ZZ) \to \GL(4,\CC)$.
				
				 We are left with determining the image of $\varphi$. One may check that $(ABDE)^3 = -I$, so that $\ker \varphi$ does not contain $-I$. By the congruence subgroup property for $\Sp(4,\ZZ)$, we conclude that $\varphi$ is injective. Since $\rho \otimes \chi$ is conjugate to the standard representation $\std$ of $\Sp(4,\ZZ)$, we conclude that $\rho$ is conjugate to a character twist of $\std$.
			\end{proof}
		\subsubsection{$g \geq 3$}
		We first establish that the image of a nonseparating Dehn twist under a low rank representation of the mapping class group is unipotent. The following lemma is due to Korkmaz, but we include the proof for completeness. We will later generalize it to representations of higher rank. 
			
		\begin{lemma}\cite[Step 1 on p.930]{korkmaz}\label{taunipotent}
			Let $\rho: \Mod(g, n) \to \GL(2g, \CC)$ be a representation with $g \geq 3$ and $n \geq 0$. Let $a$ be a nonseparating simple closed curve. Then $\rho(T_a)$ is unipotent, and if $\rho(T_a) \neq I$, then the codimension of the $1-$eigenspace is 1. 
		\end{lemma}
			The proof relies on the following result of Korkmaz, which is a consequence of the lantern relation. 
			 \begin{lemma}\cite[Lemma 8.2]{korkmaz}\label{usefulcalc6}
			 	Let $g \geq 3, n \geq 0$ and let $\rho: \Mod(g,n) \to \GL(r,\CC)$ be a representation. Let $a$ be a nonseparating simple closed curve and suppose $\rho(T_a)$ has one eigenvalue $\lambda$ whose eigenspace has codimension one. Then $\lambda = 1$. 
			 \end{lemma}
			
			\begin{proof}[Proof of \Cref{taunipotent}]
				Let $R$ be the complement of a regular neighborhood of $a$. If $A:= \rho(T_a)$ has at least three eigenvalues, then $\Mod(R)$ preserves the generalized eigenspaces of $A$, hence $\Mod(R)$ preserves a flag of subspaces each with dimension at most $2(g-1)-1$, hence $\rho|_{\Mod(R)}$ is trivial (when $g-1 = 2$, we may instead conclude by noting that $\rho([\Mod(R),\Mod(R)]) = I$, and the normal closure of the commutator subgroup of $\Mod(R)$ in $\Mod(3,n)$ is $\Mod(3,n)$). Since $\Mod(R)$ contains a Dehn twist that is conjugate to $T_a$ we conclude $A$ is trivial. 
				
				If $A$ has two eigenvalues, let $\lambda$ be the eigenvalue with greater multiplicity. First consider when $g \geq 4$. The same argument as in the previous paragraph shows that $\rho|_{\Mod(R)}$ is trivial unless the dimension of the $\lambda$ generalized eigenspaces is $2g-1$. In that case, the generalized eigenspace of the remaining eigenvalue is $1$-dimensional, hence that remaining eigenvalue must be $1$. On the other hand, \Cref{usefulcalc6} shows that if the $\lambda$-eigenspace is codimension one, then $\lambda = 1$.

				When $g =3$ and $A$ has two distinct eigenvalues, the eigenvalue $\lambda$ of greater multiplicity has multiplicity four or five. First consider when $\lambda$ has multiplicity five. Let $b$ be a nonseparating simple closed curve intersecting $a$ once. If $E_\mu^a = E_\mu^b$ we are done by \Cref{preserveespaces}, so the $\mu$-eigenspaces must be different. Similarly, if $c$ is a nonseparating simple closed curve disjoint from $a$, then $\rho(T_c)$ acts on $E_\mu^a$ either by $\mu$ or by $\lambda$. If by $\mu$, since this is true for any nonseparating simple closed curve disjoint from $a$, and in particular holds for a pair of such curves intersecting a t a point, we have that $E_\mu^a$ is $\Mod(3,n)$-invariant, and we are done, so we must have $E_{\mu}^c \subset (E_\lambda^a)_{gen}.$ Let $T$ be the complement of a regular neighborhood of $a \cup b$. Then $\Mod(T)$ preserves $W =(E_\lambda^a)_{gen} \cap (E_\lambda^b)_{gen}$, which as dimension four. Since $\dim V/W = 2$, $\Mod(T)$ acts on it by a nontrivial character. On the other hand, $\Mod(T)$ acts on $W$ either by a character or by a twist of the standard symplectic representation by \Cref{sp4rigid}. In the former case, $\rho$ annihilates the commutator subgroup of $\Mod(T),$ whose normal closure in $\Mod(3,n)$ is all of $\Mod(3,n)$. In the latter case, $T_c$ acts on $W$ by $\lambda$, which contradicts the previous observation that $T_c$ acts on $(E_\lambda^a)_{gen}$ with two distinct eigenvalues. 
				
				We are left to consider when $\lambda$ has multiplicity four and $\mu$ multiplicity two. Since $\Mod(R)$ preserves the generalized eigenspaces of $A$, we know $\Mod(R)$ acts on $(E_\mu^a)_{gen}$ with cyclic image, and on $(E_\lambda^a)_{gen}$ either with cyclic image or by a twist of the symplectic representation. In the former case, $\rho$ annihilates the commutator subgroup of $\Mod(R)$, hence $\rho$ is trivial. In the latter case, we must have $E_\lambda^a \cap E_\lambda^b \neq \emptyset$ since this holds for a pair of simple closed curves in $R$ intersecting at a point, but on the other hand, this produces a contradiction because $\Mod(R)$ preserves the filtration $E_\lambda^a \cap E_\lambda^bb \subset (E_\lambda^a)_{gen}$, which has graded pieces with dimensions 1,3, and 2. 
			\end{proof}
			
			We may now reprove the classifications of representations of the mapping class group of rank at most $2g$. 
			\begin{theorem}\label{mgagtorelli}
				Let $\rho: \Mod(g, n) \to \GL(2g, \CC)$ be a representation with $g \geq 3$ and $n \geq 0$. Then either $\rho$ is trivial, or $\rho$ is conjugate to the standard representation of $\Sp(2g,\ZZ)$. 
			\end{theorem}
			\begin{proof}
				Let $a$ be a nonseparating simple closed curve, let $b$ be a nonseparating simple closed curve intersecting $a$ once, and let $R$ be the complement of $a \cup b$. Let $a= e_1, e_2, \dots, e_n$ be simple closed curve homologous to $a$ such that $a_i$ and $a_{i+1}$ bound the $i$th marked point
				
				By \Cref{taunipotent}, $A:= \rho(T_a)$ is unipotent, and if $\rho$ is nontrivial, its $1$-eigenspace has codimension one. We now deform $\rho$ so that it underlies a $\CC$-variation of Hodge structures. 
				
				Consider the divisor $\Delta_0$ of irreducible nodal curves in $\bar{\shm}_{g,n}$. Let $N^\times$ be a punctured tubular neighborhood of $\Delta_0$. This is an $S^1$ bundle over $\Delta_0$, with fiber generated by the monodromy around $\Delta_0$, which we may take to be $T_a$. We have a short exact sequence
				$$1 \to \langle T_a \rangle \to \pi_1(N^\times) \to \pi_1(\Delta_0) \to 1.$$
				
				Since $\rho$ underlies a $\CC$-variation of Hodge structures and $A$ is unipotent, a limit mixed Hodge structures over $\Delta_0$ exists and $\rho(\pi_1(N^\times))$ is contained in a parabolic subgroup $P$ of $\GL(2g, \CC)$ which stabilizes the monodromy weight filtration $W_\bullet$. There is an induced map $\psi: \pi_1(\Delta_0) \to L$, where $L$ is the Levi subgroup associated to $P$. Since the Weight filtration satisfies the Hard Lefschetz Theorem, the largest rank of a factor of $L$ is at most$2g-2$. 
				
				If $g = 3$, $\psi$ is either a quotient of $\ZZ/10\ZZ$, or $L = \Sp(4) \times \GG_m \times \GG_m$ and $\psi$ is conjugate to one of $\std$ or $\std \otimes \lambda$ where $\lambda^{10} = 1$. If $g \geq 4$, $\psi$ is either trivial or $L = \Sp(2g-2) \times \GG_m \times \GG_m$ and $\psi$ is conjugate to $\std$. These cases occur, respectively, if the rank of the nilpotent part of $\rho(T_a)$ is greater than one or equal to one. Henceforth we assume $\rk \rho(T_a) -I = 1$.
				
				We identify $\pi_1(\Delta_0) \simeq \Mod(g-1, n+2)/ \ZZ/2\ZZ$ as the mapping class group of a genus $g-1$ surface $R$ with $n+2$ marked points, obtained by collapsing the curve $a$, normalizing the resulting surface, and forgetting the labeling of the two points which were pinched to a node. In the cases where $\psi$ has finite image, we may lift a Dehn twist around a nonseparating simple closed curve $d \subset \Sigma_{g,n} \setminus N(a)$ to $\pi_1(N^\times)$ to see that $\rho(T_d)$ is trivial. 
				
				If $\psi$ is conjugate to $\std$ (or $\std \otimes \sgn$ when $g=2$) then we may consider a homology basis $a_1, b_1, \dots, a_{g-1}, b_{g-1}$ of $\Sigma_{g-1,n+2}$ and lift the matrices $\psi(T_{a_1}), \dots \psi(T_{b_{g-1}})$ to $\pi_1(N^\times)$ and write them as
				$$\rho(T_{a_i})=
				\begin{bmatrix}
					1 & 0 &0\\
					0 & \psi(T_{a_i}) & 0\\
					0 & 0 &1
				\end{bmatrix},
				\qquad
				\rho(T_{b_i})=
				\begin{bmatrix}
					1 & 0  & 0\\
					0 & \psi(T_{b_i})\\
					0 & 0 & 1\\
				\end{bmatrix}$$ 
				
				Note that there are no off-diagonal terms. A priori, these matrices would be upper triangular, but they all commute with $A$ and $B$. In the basis respecting the monodromy weight filtration, we have $$A = \begin{bmatrix}
					1 & 0 & 1\\
					0 & I_{2g-2} & 0\\
					0 & 0 & 1\\
				\end{bmatrix}.$$
				
				By \Cref{usefulcalcs} , we may then write 
				$$B=
					\begin{bmatrix}
					1 & 0 & 0\\
					0 & I_{2g-2} & 0\\
					-1 & 0 & 1\\ 
				\end{bmatrix}.$$ Here we have used that $\img B-I \simeq \CC^{2g}/\ker(A-I)$. By the same calculation, we have $\rho(T_{e_1}) = \dots = \rho(T_{e_{n}}) = A$. From there the remainder of the proof of \Cref{sp4rigid} applies verbatim and we conclude that $\rho$ is conjugate to $\std$ (note that for $g =3, \Mod(3,n)$ is perfect, so it cannot map to a character twist of the standard representation). 
				
				We have so far assumed that $\rho$ underlies a variation of Hodge structures, but by \Cref{rigidadrho} the standard representation is strongly cohomologically rigid, so the deformation to a $\CC$-VHS was trivial to begin with. Similarly all of the representations with finite image factor through the abelianization and are rigid. 
			\end{proof}
			
\section{Uniqueness of the Torelli Map}
	In this largely expository section, we recall how a classificaiton of low dimensional representations of $\Mod(g,n)$ is enough to classify variations of Hodge structure on $\shm_{g,n}$. Throughout we will consider $g \geq 2$. 
	
	First, a word on categories. We treat $\shm_{g,n} = [\sht_{g,n}/\Mod(g,n)]$ and $\sha_h = [\HH_h/\Sp(2h,\ZZ)]$ as analytic stacks. While the general definition of a morphism of Deligne-Mumford stacks, is more complicated, for our purposes, since both spaces are quotients of contractible domains with finite stabilizers, a morphism $f: \shm_{g,n} \to \sha_h$ can be taked to be a map $\tilde{f}: \sht_{g,n} \to \HH_g$, together with a representation $\rho: \Mod(g,n) \to \Sp(2h,\ZZ)$ such that $$\tilde{f}(\gamma \cdot x) = \rho(\gamma) \tilde{f}(x).$$ 
	
	\subsection{Finite Image Representations Give Constant Maps}
		We first consider the representations $\rho: \Mod(g,n) \to \GL(r,\CC)$ with finite image. Recall that the case of a nontrivial finite image only occurs when $g=2$; for $g\geq 3$ such representations have trivial image. 
		
		By \Cref{2g-1trivial}, if $f: \shm_{g,n} \to \sha_h$ with $h<g$, we know that the underlying monodromy representation is trivial (or has finite image when $g=2$). When $g \geq 3$, the monodromy representation is trivial, so $f$ lifts to a holomorphic map $\tilde{f}:\shm_{g,n} \to \HH_h$; since $\HH_h$ is a bounded symmetric domain it follows that $\tilde{f}$ is constant. When $g=2$, $\Mod(2)$ can act via a combination of $\pm 1$ and a primitive $5$th root of unity. If the action is just by $\pm 1$, then $f$ lifts to $\HH_h$, since $-I$ acts trivially on $\HH_h$. For the other characters, after taking a degree five cover $\shm'$ of $\shm_{2,n}$, we get a constant map $\shm' \to \HH_h$. It follows that $f$ itself is constant on points, and we have proved:
		
		\begin{proposition}\label{finiteimmeanstrivial}
			Let $f: \shm_{g,n} \to \sha_h$ be a holomorphic map with $h < g$, $g\geq 3$. Then $f$ is constant. If $g = 2$, $f$ is constant on points, but can act by nontrivial automorphisms coming from $H_1(\Mod(2),\ZZ)$.
		\end{proposition}
	
		When $h = g$, the same argument applies to maps $\shm_{g,n} \to \sha_h$ whose underlying representation has finite image. 
		\begin{remark}
			It is interesting to note that When $g=2$, the finite-image representation $\rho: \Mod(2) \to \Sp(4,\ZZ)$ imposes some restrictions on the image. The reason is that $f$ respects the automorphism groups of the image. The $\rho$-equivariance of the period map forces $\img \rho$ to stabilize whatever Hodge structure is its image. If $\img \rho = \pm 1$, this is no issue, as every Abelian variety contains the $-1$ involution. But if the image is a 5th or 10th root of unity, this is very rigid. We have that $\shm_2^{ct} \simeq \sha_2$ as stacks, so the image of $f$ lies either in the Jacobian locus or in the decomposable locus of $\sha_2$. No elliptic curve over $\CC$ has an automorphism of order 5, so already the image of $f$ must lie in the Jacobian locus. Theorem 10 of \cite{onlyonesurface} states that over $\CC$ there is a unique genus two curve whose automorphism group has an element of order 5, given by the equation $y^2 = x^5-1$.
		\end{remark}
	\subsection{Rank 2g With Infinite Image is Torelli}
		We now consider maps $f: \shm_{g,n} \to \sha_g$ whose underlying representation has infinite image. Let $J: \shm_{g,n} \to \sha_g$ denote the Jacobian map. We will first give a quick conclusion for the rigidity of the Torelli map using Deligne's semisimplicity theorem. This argument is a mild simplifiaction of Hain's argument. After, we will describe Farb's argument. Both arguments, at their core, use the harmonicity of period maps and the uniqueness of a harmonic representative within a homotopy class of maps to a period domain. 
		
		We are taking as a given the conclusion to \Cref{mgagtorelli}, which states that any representation $\Mod(g,n) \to \GL(2g,\CC)$ with infinite image is conjugate to the standard symplectic representation on homology, potentially twisted by a nontrivial character when $g=2$. 
		\begin{theorem}\label{globalrigidity}
			Let $f: \shm_{g,n} \to \sha_g$ be a nonconstant morphism. When $g \geq 3$, $f = J$. If $g = 2$, $f = J$ or $J^{-}$, where the latter is the map which on points sends a curve to its Jacobian, and whose underlying local system is $\std \otimes \sgn$. 
		\end{theorem}
			\begin{proof}
				Let $\VV_f$ and $\VV_{J}$ be the local systems underlying the maps $f$ and $J$. Both underlie polarizable $\ZZ-$variations of Hodge structures. If $g \geq 3$, \Cref{mgagtorelli} implies $\VV_f$ and $\VV_J$ are isomorphic as $\ZZ$-local systems, with the isomorphism provided by the flat section of $H^0(\shm_{g,n}, \Hom_{LocSys}(\VV_f, \VV_J)) \simeq \ZZ$. By \Cref{Delignerigidity}, $\VV_f$ and $\VV_J$ each supports a unique variation of Hodge structures. Therefore, they induce a unique weight $0$ variation of Hodge structures on $\Hom_{LocSys}(\VV_f,\VV_J)$, so that the local system isomorphism is in fact an isomorphism in the category of variations of Hodge structures. 
				
				By working on moduli spaces with level structures, we may upgrade the isomorphism $f\simeq g$ to an equality. For $\ell \geq 3$, $\sha_g[\ell] = [\HH_g/\Sp(2g,\ZZ)(\ell)]$ is a fine moduli space. Moreover, working on the cover $\shm_{g,n}[\ell] = [\sht_{g,n}/\Mod(g,n)[\ell]]$, where $\Mod(g,n)[\ell] = \ker(\Mod(g,n) \to \Sp(2g, \ZZ/\ell\ZZ))$, $f$ and $g$ lift to maps $\tilde{f}, \tilde{g}: \shm_{g,n}[\ell] \to \sha_{g}[\ell]$. The argument from the previous paragraph yields $\tilde{f} \simeq \tilde{g}$, which, since the moduli spaces are fine, implies $\tilde{f} = \tilde{g}$. This also means that $f$ and $g$, as maps to the coarse space $A_g$, are equal since $\tilde{f}$ and $\tilde{g}$ agree on points.
				
				When $g = 2$, the symplectic local system $\VV_f$ can be a twist of $\VV_J$ by the sign character. To get around that, we can twist $\VV_f$ so that the image of a nonseparating Dehn twist is unipotent, and we handle the other case by considering $\VV_f \otimes \chi$. We have arranged so that the same argument as when $g \geq 3$ gives $\VV_f \simeq \VV_J$ as Hodge structures since $\VV_f$ is now the local system which is conjugate to $\VV_J$, but now there is no reason  $\Hom_{LocSys}(\VV_f \otimes \chi, \VV_W))$ should be nonzero. In fact, using $\VV_f \simeq \VV_J$, we have $\Hom_{\CC}(\VV_f \otimes \chi, \VV_W) \simeq \End_\CC (\VV_J) \otimes \chi^\vee$. A flat section therefore is a morphism of local systems $\chi \to \End_\CC(\VV_J)$. Since $\VV_J = \Sym^2 \VV_J \oplus \wedge^2 \VV_J \oplus 1$ as a representation of $\Sp(4,\ZZ)$ and $\chi$ does not occur in this decomposition, there is no such morphism. It follows that on fine moduli spaces, $f$ and $J$ are not equal if $f_* \simeq \rho \otimes \chi$. 
				
				At the level of stacks, we may conclude $f$ is isomorphic to one of $J$ or $J^{-}$, and the same argument with level structures will give either $\tilde{f} = \tilde{J}$ or $\tilde{f} = \tilde{J}^-$ as maps $\shm_{2,n}[\ell] \to \sha_g[\ell]$ with $\ell \geq 3$. On the other hand, since the character $\sgn$ acts trivially on $\HH_2$, $f, J$, and $J^{-}$ all descend to the same map to the coarse space $A_2$.
			\end{proof}

		We now describe how Farb conclude that $f=J$ as maps to $A_g$. We treat the case $g \geq 3$ first, to which Farb's argument in \cite{farbrigid} applies verbatim, so we only sketch it briefly. By the classification from \Cref{mgagtorelli}, any representation $\Mod(g,n) \to \GL(2g,\CC)$ is conjugate to the representation given by the Torelli map. Since there exist complete curves in $\shm_{g,n}$ whose fundamental groups surject onto $\Mod(g,n)$, it is enough to restrict the homotopy to such a curve: $H: C \times [0,1] \to \sha_g$. Specifically, restrict $H$ to a (necessarily non-complete) Lefschetz curve $C$. The monodromy group of $\rho|C$ is simple (the identity component of its its algebraic closure is $\Sp(2g)$), and the monodromy around the punctures of $C$ is given by a Dehn twist about the vanishing cycle, hence it has infinite image. The work of Saito then implies that for any $t$, $H_t: C \to \sha_g$ is an isolated point in $\Hom(C, \sha_g)$, so the homotopy must be constant. Hence $f = J$. 
		
		When $g=2$, the classification of representations of $\Mod(2,n)$ with infinite image also allowed for character twists of the standard representation. Starting with a representation $\rho: \Mod(2,n) \to Sp(4,\ZZ)$, this means we can twist $\rho$ by the sign character, since the image of a nonseparating Dehn twist has a unique eigenvalue, so up to conjugacy, there are two representations $\Mod(2,n) \to \Sp(4,\ZZ)$. As noted above, the involution $\chi: \Mod(2,n) \to \{\pm I_4\}$ acts trivially on $\HH_2$, so even if $f_* = \std \otimes \sgn$, the lift $\tilde{f}: \sht_{2,n} \to \HH_2$ is the same whether or not there is a twist; in particular, the lift is the same as the lift of the Torelli map. While $\shm_2$ does not have any complete curves (it is affine), any Lefschetz curve in $\shm_2$ with fundamental group surjecting onto $\Mod(2)$ will have a simple monodromy group. Farb's argument therefore gives that either $f = J$ or $f = J^-$, where $J^-$ is the same map on points, but acts via $-I$ on automorphisms. In particular, the induce the same map to the coarse space $A_g$. 
	
\section{Maps $\shm_{g,n} \to \sha_{2g-1}$}
	This section begins the more algebro-geometric portion of the paper. The goal is still to classify $\shm_{g,n} \to \sha_h$, but instead of proceeding by classifying all representations of the mapping class group, we begin with the morphism of stacks and use the constraints of such a map being a variation of Hodge structure. As with the similar classification theorems of Korkmaz \Cite{korkmaz} and Kaufmann-Salter-Zhang-Zhong \Cite{kaufmann2026linearrepresentationsmappingclass}, we often will argue by restricting our representation to a subsurface. However, the Hodge-theoretic context allows us to consider a related representation, which is the limit mixed Hodge structure around the boundary divisor $\Delta_0$, which will be a representation of the mapping class group $\Mod(g-1, n+2)$. This representation has the virtue of still being a variation of Hodge structures, unlike the general restriction to a subsurface.
	
	Before analyzing the degeneration of the variation of Hodge structures, we must do some preliminary work on classifying the eigenvalues of the local monodromy around the divisor $\Delta_0$ (equivalently, the image of a Dehn twist around a nonseparating curve.) For most cases under consideration, this is already known to be unipotent by \Cite[Proposition 5.1]{kaufmann2026linearrepresentationsmappingclass}. Specifically, they show any representation $\Mod(g,n) \to \GL(r,\CC)$ with $g\geq 4$ and $r \leq 4g-3$ sends a Dehn twist about a nonseparating simple closed curve to a unipotent matrix. For our theorem, then, we must prove the same result for $g =3$ and when $r=4g-2$; additionally we must handle the fact that when $g=2$ the image of a nonseparating Dehn twist does not have to be unipotent. 
	
	The argument is long and some lemmas handle cases for multiple parts of the argument, so we provide the following outline:
	
	\begin{itemize}
		\item Restrict the possible eigenvalues. When $g= 2$, \Cref{g2oneeval} shows the local monodromy around $\Delta_0$ has either $1$ or $-1$ as its unique eigenvalue in our range. For $g \geq 3$, \ref{mgagunipotent} shows that the local monodromy around $\Delta_0$ is unipotent in our range.
		\item Base case when $g=2$. \Cref{g2agunique} gives another proof of the uniqueness of the Torelli map $\shm_2 \to \sha_2$. 
		\item Getting rid of points. \Cref{nork1monodromy} reduces the classification maps $\shm_{2,n} \to \sha_2$ to the case where $n=0$, completing the classification for $g=h=2$. Combined with \Cref{nork2monodromy}, it also reduces the classification of maps $\shm_{2,n}\to \sha_3$ to the case where $n=0$.
		\item Maps $\shm_2 \to \sha_3$. There are 16 cases to consider, depending on the rank of the local monodromy around the divisors $\Delta_0$ and $\Delta_1$ each having ranks in $\{0,1,2,3\}$. Some cases are dispensed with in \Cref{linalgm2a3} by finding an invariant vector and using the theorem of the fixed part. Other cases we rule out in \Cref{biratm2a3}, showing that the pullback of the Hodge class on $\sha_3$ is trivial, hence the map is trivial. The only case left is when the local monodromy around $\Delta_0$ has nilpotent rank 1 and the local monodromy around $\Delta_0$ is trivial, for which $f$ factors through the Torelli map $\shm_{2,n} \to \sha_2$. By the theorem of the fixed part and highest weight considerations, $f$ has the desired factorization.
		\item Inducting. Having completed the base case, we proceed with the induction on $(g,h)$ in lexicographic order (always with $h \leq 2g-1$) and $n$ arbitrary.
	\end{itemize}
	
	In every case, we will show that the underlying representation is conjugate to a direct sum of the standard symplectic representation (or a twist thereof when $g=2$) with a trivial representation (or a finite cyclic one when $g=2$). By the argument of \Cref{globalrigidity}, this is sufficient to show that the varying part of $f$ equals the Torelli map on points, or, on a level cover, an equality of weight on $\ZZ$-variations of Hodge structures. There is, of course, freedom in the constant Abelian variety attached to $f$.

	\subsection{Eigenvalues of Families of Abelian Varieties over $\shm_{g,n}$}
		Given a family of abelian varieties $f: \shm_{g,n} \to \sha_h$, there is a corresponding representation of the mapping class group $\rho:\Mod(g,n) \to \Sp(2h,\ZZ)$. By definition this underlies an integral polarized weight one variation of Hodge structures, so the image of a nonseparating Dehn twist $T_c$ has a Jordan decomposition $\rho(T_c) = S(I+N)$ with $S^k = I$ and $N^2 = 0$. From $N^2 = 0$, one immediately sees that $\rk N \leq h$, and the eigenvalues of $\rho(T_c)$ are controlled by partitions of $2h$ compatible with the rank of $N$. Moreover, for a symplectic representation, eigenvalues come in reciprocal pairs. We first show that when $g\geq 3$ and $h \leq 2g-1$, $\rho(T_c)$ is forced to be unipotent. Meanwhile when $g = 2$ eigenvalues of $\pm 1$ are possible.
		
		\begin{lemma}\label{g2oneeval}
			Let $f: \shm_{2,n} \to \sha_h$ be a nonconstant morphism with $h = 2,3$. Then the image of a nonseparating Dehn twist under the representation of $\Mod(2,n)$ associated to $f$ has one eigenvalue, equal to either $1$ or $-1$.
		\end{lemma}
			\begin{proof}
				The case $h=2$ follows from \Cref{sp4rigid} after observing that the semisimple part of $\rho(T_c),$ $\lambda I$, lies in the center of $\Sp(4,\ZZ)$, hence $\lambda = \pm 1$. Alternatively, to avoid logical dependence on that result, first note that if $\lambda$ is an eigenvalue of $\rho(T_c)$ with multiplicity one, then, if $\lambda \neq \pm 1,$ $\lambda\inv$ is also an eigenvalue. In particular, in order for $\rho(T_c)$ to have an eigenvalue $\lambda \neq 1$ and not be semisimple, $\lambda$ must have multiplicity two. On the other hand if $\lambda = \pm 1$, the $1$ and $-1$ eigenspaces are orthogonal and therefore must have the same dimension. Fixing $\rho(T_c)$ in its Jordan basis, let $a$ and $b$ be two nonseparating simple closed curves intersecting once, each disjoint from $c$. 
				
				If $\lambda \neq \pm 1$ and there is more than one eigenvalue, the one can have spectra $(\lambda \lambda, \lambda\inv, \lambda\inv)$ or $(\pm 1, \pm 1, \mp 1, \mp 1)$, while if $h = 3$, the possible spectra are $(\lambda, \lambda, \lambda, \lambda\inv, \lambda\inv, \lambda\inv)$, $(\pm 1,\pm 1,\lambda,\lambda,\lambda\inv, \lambda\inv)$, $(1,1,-1,-1,-1,-1)$, or $(-1,-1,1,1,1,1)$. 
				
				We now show all of these cases are impossible by what is essentially the argument of \Cref{g2finiteimage}. With the same curves $a$ and $b$ as before, $\rho(T_a)$ and $\rho(T_b)$ are block diagonal with respect to the Jordan decomposition of $\rho(T_c)$. Then the braid relation between $\rho(T_a)$ and $\rho(T_b)$ implies the two matrices are equal on any Jordan block of size $2$ or $3$, so that those generalized eigenspaces are $\Mod (2,n)$-invariant.  
				
				We have thus shown in all instances that $\rho(T_c)$ has eigenvalues other than $\pm 1$ that $\rho$ is reducible into representations of rank at most $3$, which all necessarily have finite image. 
			\end{proof}
			
		We now prove a similar result for $g=3$. Much of the argument of \Cref{g2oneeval} applies verbatim. For the calculations in the following lemma, we will need to explicitly construct such filtrations. Once we know that the image of a nonseparating Dehn twist is unipotent, we will reinterpret this lemma as a statement about the limiting mixed Hodge structure associated to a period map on $\shm_{g,n}$. 
		
		\begin{lemma}
			Let $f: \shm_{3,n} \to \sha_h$ be a nonconstant morphism with $3 \leq h \leq 5$. Let $\rho$ be the corresponding representation of $\Mod(3,n)$ and let $c$ be a nonseparating simple closed curve. Then $\rho(T_c)$ is unipotent. 
		\end{lemma}
			\begin{proof}
				We will consider each value of $h$ individually, taking the classification of representations $\Mod(g,n) \to GL(r,\CC)$ for $r \leq 2g-1$ and $g=2,3$ for granted. Let $c$ and $d$ be two nonseparating simple closed curves meeting at one point, and let $R$ be a the complement of a regular neighborhood of $a \cup b$. $\Mod(R)$ commutes with $\rho(T_c)$, so it preserves the generalized eigenspaces of $\rho(T_c)$. By \Cref{usefulcalcs}, if a Jordan block has size 1 or 2, the braid relation implies that block is $\Mod(3,n)$-invariant, and we are done by \Cref{2g-1trivial}.
				
				For the rest of the analysis, recall that eigenvalues of symplecitc matrices come in reciprocal pairs. Moreover, since $\rho$ comes from a variation of Hodge structures, the eigenvalues are also roots of unity. If $p(t)$ is the characteristic polynomial of $\rho(T_c)$, then $p(t)$ must be a product of cyclotomic polynomials. Since the degree of a cyclotomic polynomial is $\varphi(n)$ where $n$ is the root of unity and $\varphi$ is the totient function, and since we may assume a generalized eigenspace for an eigenvalue has dimension at least $3$, we see that $3\varphi(n) \leq 2h$. For $h \in \{3,4,5\}$, this implies $\lambda$ is a primitive $n$th root of unity for $n \in \{1,2,3,4,6\}$. 
				
				Now suppose $\lambda \neq 1$ is an eigenvalue of $\rho(T_c)$. Let $d$ be the dimension of the generalized $\lambda-$eigenspace. If $d = 3$, consider $W = E_{\lambda}^{gen} \oplus E_{\lambda\inv}^{gen}$. Note that $\lambda \neq -1$ since its generalized eigenspace is odd dimensional. The representation $\Mod(R) \to \Sp(W)$ is reducible and has finite cyclic image. For $h = 3,4$, we are done, as $\rho|_{\Mod(R)}$ annihilates the commutator subgroup of $\Mod(R)$. When $h = 5$, $V/W = E_1^{gen}$ has dimension four and $\Mod(R)$ can act on it as the standard symplectic representation. If so, we conclude that $\dim E_1^c = 3$, hence $\Mod(R)$ preserves the flag $0 \subset E_1^c \subset (E_1^c)_{gen} \subset (E_\lambda^c)_{gen} \subset (E_{\lambda\inv}^c)_{gen}$, whose graded pieces all have dimension at most 3, so $\rho|_{\Mod(R)}$ annihilates the commutator subgroup of $\Mod(R)$, so that $\rho$ was actually trivial. 
				
				If the generalized $\lambda-$eigenspace has dimension 4, we have a representation $\phi: \Mod(R) \to \GL(E_\lambda^{gen})$. If this representation is trivial then $E_\lambda$ is $\Mod(3,n)$ invariant, and we are done by induction. If nontrivial, then $\lambda = -1$, and $\phi$ is known to be a scalar twist of the symplecitc representation.  $E_\lambda^c \cap E_\lambda^d \subset E_\lambda^c \subset (E_\lambda^c)_{gen} \oplus (E_{\lambda\inv}^c \cap E_{\lambda\inv}^d) \subset (E_\lambda^c)_{gen} \oplus (E_{\lambda\inv}^c)_{gen} $ is a filtration with graded pieces having ranks at most $3$, so we are done by \Cref{flag}. 
				
				If the generalzied $\lambda$-eigenspace has dimension $5$, we have a representation $\Mod(R) \to \GL(E_\lambda^{gen})$. Taking determinants, we must have $\lambda^5 = 1$, hence $\lambda = 1$
				
				This concludes the analysis. An eigenvalue $\lambda \neq 1$ has a generalized eigenspace of dimension at least 3, while a generalized $-1$ eigenspace must have even dimension, hence at least $4$. Hence the only nontrivial cases, where $\rho$ does not immediately have an invariant subspace, are when a generalized eigenspace has dimension $h$, which we handled above. 
				\end{proof}
				
		We can finally prove the result for higher $g$ and $h$, where less case-by-case analysis is required. We adapt the argument of Kaufmann et. al. (\Cite[Proposition 5.1]{kaufmann2026linearrepresentationsmappingclass}), using the additional assumption that our representation comes from a family of abelian varieties. 
	
		\begin{lemma}\label{mgagunipotent}
			Let $f: \shm_{g,n} \to \sha_h$ be a nonconstant morphism with $g \geq 3$ and $g \leq h \leq 2g-1$. Let $\rho$ be the corresponding representation of $\Mod(g,n)$ and let $c$ be a nonseparating simple closed curve. Then $\rho(T_c)$ is unipotent. 
		\end{lemma}
		\begin{proof}
			We may assume $\rho(T_c)$ is not of the form $\lambda I$ with $\lambda \neq 1$ because, if that were the case, we may find seven disjoint nonseparating simple closed curves satisfying the lantern relation; applying $\rho$ to the lantern relation yields $\lambda^4 = \lambda^3$. We will induct on $h$ and $g$, with $g=3$ and $h=g$ as the base case, treated by the previous lemma. First let $\lambda \neq \pm 1$ be an eigenvalue. The arguments of the previous lemmas apply, so that the (generalized) $\lambda-$eigenspace has dimension at least $3$, otherwise the inductive hypothesis applies. Denote by $E_\lambda^{(gen)}$ the (generalized) $\lambda$-eigenspace of $\rho(T_c)$. Since $\lambda \neq \pm 1$, $\lambda\inv$ is also an eigenvalue, and its eigenspace has the same dimension. Thus $3 \leq \dim E_\lambda \leq h$. If there is a third eigenvalue in addition to $\lambda$ and $\lambda\inv$, then in fact $\dim E_\lambda \leq h-2 \leq 2g-3$, at which point \Cref{2g-1trivial} applies: if $R$ is the complement of an annular neighborhood of $c$, then $\Mod(R) \to \GL(E_\lambda)$ is the trivial representation, contradicting that $\lambda \neq 1$. 
			
			Thus if $\lambda \neq \pm 1$ is an eigenvalue, we only have to consider when $\dim E_\lambda^{gen} = h$ and $\rho(T_c)$ only has eigenvalues $\lambda$ and $\lambda\inv$. In this case, $\Mod(R) \to \GL(E_\lambda^{gen})$ is a $2g-1$ dimensional representation of $\Mod(R)$ with a unique eigenvalue $\lambda$. If this representation is trivial, the lantern relation implies $\lambda = 1$. Otherwise, $E_\lambda^{gen} \neq E_{\lambda}$, and the representation $\Mod(R) \to \GL(E_\lambda)$ is known by induction, and in particular satisfies $\lambda = 1$.
			
			We finally consider when $\lambda = -1$. We have $4 \leq \dim E_{-1} \leq 2h-4$. In particular if $\dim E_{-1} > 2g-2$ (hence $\geq 2g$, since it has even dimension), we have $\dim \CC^{2h}/E_{-1} \leq 2h-2g$. Since $h \leq 2g-1$, we conclude either $\dim E_{-1}$ or $\dim \CC^{2g}/E_{-1}$ is at most $2g-2$-dimensional. These spaces are invariant under $\Mod(R)$, where we see they have a unique eigenvalue equal to $1$, a contradiction. 
		\end{proof}
	\subsection{The Base Case $g=2$, Again}
		We now give another argument for the uniqueness of the Torelli map in genus $2$. We now start with a holomorphic map $f:\shm_2 \to \sha_2$ and show that $f$ is the Torelli map. By \Cref{g2oneeval}, we may twist the representation $f_*$ by the sign character to assume the local monodromy around the boundary of $\bar{\shm_2}$ is unipotent, so that $f$ extends to the Deligne-Mumford and Satake compactifications without taking any cover. 
		
		The local monodromy around $\Delta_0$ and $\Delta_1$ are given by Dehn twists about nonseparating and separating curves, respectively. If $a$ and $b$ are two nonseparating simple closed curves meeting at a point, a separating twist $T_\delta$ is given by $(T_aT_b)^6$. Writing $f_*(T_a) = I+a$ and $f_*(T_b) = I+b$, with $a^2 = b^2 =0$, one has $f_*(T_\delta) = I+12(a+aba)$. For $\shm_{g,n}$, the same holds for the local monodromy around $\Delta_{1,\{0,\dots, n\}}$, parameterizing pointed curves with an elliptic tail, assuming none of the marked points are contained in the neighborhood of $a \cup b$.
		
		\begin{proposition}\label{g2agunique}
			Let $f: \shm_2 \to \sha_h$ be a map which is nonconstant on points, with $h \leq 2$. Then $h=2$ and $f$ is the Torelli map, or the Torelli map twisted by the $\pm 1$ character.
		\end{proposition}
			\begin{proof}
				By definition a map $\shm_2 \to \sha_0$ is constant. There can be no nonconstant map $\shm_2 \to \sha_1$ because $H_1(\shm_2, \CC) = 0$, so $\shm_2$ cannot have a nonconstant map to a curve.
				 By the discussion in the previous paragraph, $f$ extends to a morphism $\bar{\shm}_2 \to \sha_2^*$, where $\sha_g^*$ is the Satake compactification. We consider the resulting maps on the boundary divisors $\Delta_0$ and $\Delta_1$. If $\Delta_0$ were to map to the interior of $\sha_2$, the local monodromy around it would be trivial, but this local monodromy is generated by a nonseparating Dehn twist, so the underlying representation $f_*$ is forced to be trivial. Therefore $\Delta_0$ is forced to map to the boundary of $\sha_2$. 
				
				If the monodromy around $\Delta_1$ is trivial, then $f$ extends to $\shm_2^{ct} \simeq \sha_2$, so $f$ is an isomorphism. Otherwise, $\Delta_1$ maps to the boundary of $\sha_2$. Properness is preserved by base change, so we have a proper morphism $\shm_2 \to \sha_2^*$. This cannot occur, because $\sha_2^*$ contains complete curves, whereas $\shm_2$ is affine. 
				
				It only remains to account for the preliminary character twist if the local monodromy around $\Delta_0$ had a semisimple part of $-I$. Since $-I$ acts trivially on $\HH_2$, this does not change image of $f$ on points, so we have shown that the representation underlying $f$ is a twist of the symplectic representation by the sign character. 
			\end{proof}
	
		There are two remaining tasks in this section. First, we strengthen the prevous argument to include maps $f: \shm_2 \to \sha_3$. The remainder of the section will then be devoted to showing that maps $\shm_{2,n} \to \sha_h$ with $h = 2,3$ and $n>0$ factor through the projection to $\shm_{2}$. 
		 
		We first turn our attention to maps $\shm_2 \to \sha_3$. As the local monodromy operators around the two boundary strata may each have logarithms of rank $\{0,1,2,3\}$, there are 16 cases to consider. We cover a large number of them by linear algebra considerations, and the remaining cases we cover using the birational geometry of $\bar{\shm_2}$. 
		
		\subsubsection{Linear-Algebraic Considerations}\label{linalgm2a3}
		Let $a$ and $b$ be two nonseparating simple closed curves intersecting once. By the braid relation, a Dehn twist around a separating curve $\delta$ bounding a regular neighborhood of $ a\cup b$ is given by $T_\delta =(T_a T_b)^6$. Writing $f_*(T_a) = I+a$ and $f_*(T_b) = I+b$, the braid relation gives $f_*(T_\delta) =I+ 12(a+aba) =I+ 12(b+bab)$. On the linear algebra side, we may rule out many cases immediately. 
		
		\begin{lemma}\label{superrigidity}
			Let $f: \shm_2 \to \sha_3$ be a morphism. If the local monodromy around $\Delta_0$ is trivial, then $f$ is constant. If the local monodromy around $\Delta_0$ has nilpotent rank 1, or if the local monodromy around $\Delta_1$ is trivial, then $f$ either factors through a product locus $\{*\} \times \sha_2$ or is constant.
		\end{lemma}
			\begin{proof}
				The first statement follows because the mapping class group is normally generated by nonseparating Dehn twists. If the local monodromy around $\Delta_0$ has nilpotent rank 1, then $aba=-a$. Indeed, since $\img (a+aba) \subset \img (a)$, $a+aba$ has nilpotent rank 1 or 0. If it has rank 1, then its image line is the same image as $a$. Since $a+aba=b+bab$, the same reasoning shows that $\img a = \img b$. Since the commutator subgroup is normally generated by $T_aT_b^{-1}$, it follows that $f_*$ factors through the abelianization, so by \Cref{finiteimmeanstrivial}, $f$ is trivial. Therefore $\rk(a+aba) = 0$, so that the Torelli group lies in the kernel of $f_*$\footnote{Note that in genus $2$, the Torelli group is generated by genus one separating twists.}. $f_*$ is therefore determined by the quotient homomorphism $\Sp(4,\ZZ) \to \Sp(6,\ZZ)$. We may assume this has infinite image, so Margulis superrigidity reduces this to the classification of 6-dimensional representations of $\Sp(4, \CC)$. There are two candidates: $V_{\omega_1} \oplus 1^{\oplus 2}$ and $V_{\omega_2} \oplus 1$, where $V_{\omega_i}$ is the representation of highest weight $i$. The latter is not possible as a homomorphism to $\Sp(6,\ZZ)$ because the highest weight 2 representation preserves a symmetric form. So $f_*$ factors as the standard representation of $\Sp(4,\ZZ)$ with a trivial factor. Similarly, if the local monodromy around $\Delta_1$ is trivial, $f$ factors through the Torelli map to $\sha_2$, so the same considerations apply.
			\end{proof}
			
			We may also immediately rule out the case where the local monodromy around $\Delta_0$ has nilpotent rank $2$ while the local monodromy around $\Delta_1$ has rank 3. This is because the image of the operator $a+aba$ is contained in the image of $a$, so its rank cannot exceed the rank of $a$. Similarly if the logarithms of the monodromy operators around both divisors have rank 3, then $\rk a = \rk a+aba = \rk b = 3$. Since these are nilpotent operators on a 6-dimensional vector space, the images of $a$ and $b$ are Lagrangian subspaces. In particular $\img (a) = \ker (a)$ and $\img(b) =\ker b$, hence $ab=ba=0$. But then the braid relation $a+aba=b+bab$ implies $a=b$, so that $f_*$ kills the commutator subgroup of $\Mod(2)$, and $f$ is constant. 
			
			The last case we will rule out with linear algebra is when the local monodromy around $\Delta_0$ has nilpotent rank $2$ while the local monodromy around $\Delta_1$ has nilpotent rank 1. 
			\begin{lemma}\label{selfadjtrick}
				Let $f: \shm_2 \to \sha_3$ be a morphism where the local monodromy around $\Delta_0$ has nilpotent rank 2 and the local monodromy around $\Delta_1$ has nilpotent rank $1.$ Then $f_*$ is reducible, and consequently $f$ factors through a product locus. 
			\end{lemma}
				\begin{proof}
					We use the fact that $\Mod(2)$ is generated by five Dehn twists around nonseparating simple closed curves $a,b,c,d,$ and $e$. Assume these curves are in standard position, i.e. that $a$ and $b$ intersect once, $b$ and $c$ intersect once, etc., and the curves are otherwise disjoint. Write $f_*(T_x) = I+ x$, where $x$ has rank $2$ and $x^2 = 0$. We will show the vector $v \in \img (a+aba)$ is preserved by the images of all of these twists.
					
					 Let $z = \{c,d,e\}$. We will show $az=0$. Since both $a$ and $z$ have rank 2, the possible ranks of $az$ are 1 and 2. If $\rk az = 2$, then $\img az \subset \img a$ have the same rank, so they are equal. The same applies to $\img z,$ hence $\img a = \img z$. In particular $\img d = \img e$, so that $f_*$ kills the commutator subgroup of $\Mod(2)$. Meanwhile, if $\rk az = 1$, then we consider the symplectic adjoint $(az)^*.$ Here the adjoint of a nilpotent operator $N$ with respect to the symplectic form $\omega$ is defined by the relation $\omega(Nu,v) = \omega(u,N^* v)$. We have $a^* = -a$ and $z^* = -z$, hence $(az)^* = az$ is self-adjoint. Thus the bilinear form $b(u,v):=\omega((az)u,v)$ is alternating: $\omega((az)u,v) = \omega(u, (az)v) = -\omega((az)v,u)$, so it has even rank. But $\rk b(u,v) = \rk az = 1$, a contradiction. 
					 
					It follows that $az = 0$ for $z = \{c,d,e\}$. Since $v \in \img a$, it follows that $zv = za(x)=0$. Since $a+aba=b+bab$, we have $v= (b+bab)x$ for some nonzero $x$, so that $bv =0$ as well. Therefore $v$ is an invariant vector under $f_*$. By the theorem of the fixed part, $f_*$ has a constant subvariation. If the rank of that constant factor is more than two, $f$ is trivial, while if the rank is two, $f$ can either be trivial or factor through the surjection $\Mod(2) \to \Sp(4,\ZZ)$. 
				\end{proof}

		\subsubsection{Some Birational Geometry on $\bar{\shm}_2$}\label{biratm2a3}
			We rule out the remaining cases for maps $\shm_2 \to \sha_3$ using the birational geometry of $\shm_2$. We now switch perspectives, and think of the ranks of the local monodromy operators as determining which boundary stratum of $\sha_3^*$ the divisors $\Delta_i$ are mapped to. Indeed, if the monodromy around $\Delta_i$ has nilpotent rank $r$, then the action of $\pi_1(\Delta_i)$ on the Levi of the parabolic subgroup preserving the weight filtration induces a morphism $\Delta_i \to \sha_{3-r}$. We must consider the following cases: $(\Delta_0 \mapsto \sha_1, \sha_0;\Delta_1 \mapsto \sha_1)$ and $(\Delta_0 \mapsto \sha_0; \Delta_1 \mapsto \sha_2)$. We eliminate these possibilities in \Cref{d1noa1} and \Cref{d0noa0} below.
			
			To begin, we give a brief overview of the divisor classes on $\bar{\shm}_2$ and $\bar{\shm}_{1,n}$. Chapter 13 of \cite{ACGH2} contains most of the facts we will use.  If $\pi: \shc_g \to \shm_g$ is the universal curve and $\omega_{\pi}$ is the relative dualizing sheaf, the \textit{Hodge bundle} is $\lambda := \pi_* \omega_\pi$. The divisor class of the boundary $\bar{\shm}_{2} \setminus \shm_2$, denoted $\Delta$, decomposes into two classes, $\Delta \equiv \delta_0 + 2\delta_1$, which are divisor classes associated to the two codimension one strata of the boundary. The universal curve $\shm_{g,n} \to \shm_{g,n}$ has $n$ disjoint sections $\sigma_i$, and the pullbacks of the relative dualizing sheaf along these maps are denoted $L_i := \sigma_i^* \omega_\pi$. Their first Chern classes are called the $\psi$ classes $\psi_i := c_1(\sigma_i^* \omega_\pi)$. In $\bar{\shm}_{1,2}$, we let $\delta_{00}$ denote the class of the divisor parameterizing an elliptic curve with a 2-pointed rational tail. The only other boundary divisor in $\bar{\shm}_{1,2}$ is $\Delta_0$, whose divisor we represent by $\delta_0$, and the only boundary divisor in $\bar{\shm}_{1,1}$ is the cusp, whose class we also denote by $\delta_0$. 
			
			On $\bar{\shm}_2$, the two boundary divisors have normalizations $\tilde{\Delta_0} \simeq \overline{\shm_{1,2}}$ and $\tilde{\Delta_1} \simeq \overline{\shm_{1,1}} \times \overline{\shm_{1,1}}$. There are natural maps from these moduli spaces to $\bar{\shm}_2$, called clutching maps. For $\overline{\shm_{1,2}}$, the clutching map $\xi_0: \overline{\shm_{1,2}} \to \bar{\shm}_2$ arises by gluing the two marked points to a node, and for $\overline{\shm_{1,1}} \times \overline{\shm_{1,1}}$ the clutching map $\xi_1$ comes from gluing the two elliptic curves together at their marked points. The images of these two maps on Picard groups can be calculated combinatorially. 
			
			\begin{proposition}\label{m2barrelations}
				\begin{enumerate}
					\item The Picard group of $\bar{\shm}_2$ is spanned by the classes $\lambda, \delta_0, \delta_1$, subject to the relation $10\lambda = \delta_0 + 2\delta_1$
					\item The nef cone of $\bar{\shm}_2$ is spanned by the rays $\delta_0 + 12\delta_1$ and $\lambda$. 
					\item Let $\xi_0: \tilde{\Delta_0} \to \bar{\shm}_2$ be the clutching morphism for $\Delta_0$. Then 
					\begin{enumerate}
						\item $\xi_0^* \lambda = \lambda$
						\item $\xi_0^* \delta_0 = \delta_0 - \psi_1 - \psi_2$
						\item $\xi_0^* \delta_1 = \delta_{00}$
					\end{enumerate}
					\item Let $\xi_1: \tilde{\Delta_1} \to \bar{\shm}_2$ be the clutching morphism for $\Delta_1$. Then
						\begin{enumerate}
							\item $\xi_1^* \lambda = \lambda_1 + \lambda_2$
							\item $\xi_1^* \delta_0 = 12(\lambda_1 + \lambda_2)$
							\item $\xi_1^* \delta_1 = -(\psi_1 + \psi_2) = -(\lambda_1 + \lambda_2)$. 
						\end{enumerate}
					\item On $\bar{\shm}_{1,1}$, we have $\psi = \delta_0 = 12\lambda$ and $\deg \lambda = \frac{1}{24}$.
				\end{enumerate}
			\end{proposition}
					\begin{proof}
						The first two items are stated in \Cite{hassett} in the introduction to Section 4. The calculations there are on the coarse space $\bar{M}_{2}$, and we have normalized them to be on the stack. The next two are \Cite[Ch. 17, Lemma 4.35]{ACGH2} and \Cite[Ch. 17, Lemma 4.36]{ACGH2}. For the last statement, on $\bar{\shm}_{1,1}$ the relation $12\lambda = \delta_0$ holds (see e.g. \Cite[Ch 13, Example 7.10]{ACGH2}); the Hodge class generates $\Pic(\bar{\shm}_{1,1})$, and $\deg \delta_0 = \frac{1}{2}$, coming from the involution of the marked points on the normalization of the nodal cubic. 
					\end{proof}
			Applying the clutching formulas for the various boundary divisors and using the relations on $\bar{\shm}_{1,1}$, we have the following identities:
			\begin{corollary}
				Let $D = a\lambda + b \delta_0 + c\delta_1$. Then 
				\begin{enumerate}
					\item $\xi_0^* D = (a+10b)\lambda + (c-2b)\delta_{00}$
					\item $\xi_1^* D = (a+12b-c)(\lambda_1 + \lambda_2)$.
				\end{enumerate}
			\end{corollary}
			Throughout, we let $L$ be the ample class on $\sha_3$ and write $f^*L = a \lambda + b\delta_0 + c\delta_1$.
			
			\begin{lemma}\label{d1noa1}
				Both divisors $\Delta_0$ and $\Delta_1$ cannot map to $\sha_1$ unless $f$ is constant.
			\end{lemma}
				\begin{proof}
					Since $\dim \tilde{\Delta_1} = 2$, if $\Delta_1 \mapsto \sha_1$, there is a curve on $\tilde{\Delta_1}$ which intersects $f^*L$ trivially. Since $\Pic \overline{\shm_{1,1}} \times \overline{\shm_{1,1}} \simeq \ZZ \oplus \ZZ$ and the curve is effective, we may write it as a positive linear combination of the fiber classes $F_1 = \bar{\shm}_{1,1} \times \{*\}$ and $F_2 = \{*\}\times \bar{\shm}_{1,1}$. The intersection $F_i \cdot \lambda_j$ with $i\neq j$ is the degree of the Hodge bundle on $\bar{\shm}_{1,1}$, which is $\frac{1}{24}$. Hence we have $$0= C \cdot f^*L = (xF_1 + yF_2) \cdot \xi_1^* f^*L = (xF_1 + yF_2)(a+12b-c)(\lambda_1 + \lambda_2) = (a+12b-c)(\frac{x+y}{24}).$$ 
					
					Since $x+y \neq 0$, we must have $\xi_1^* f^*L  = 0$, so that $f$ contracts $\Delta_1$ to a point. Using the relation $a+12b-c=0$, we see that $\xi_0^* f^*L$ is a scalar multiple of $\lambda + \delta_{00}$. From \cite[17. Lemma 4.28]{ACGH2}, on $\bar{\shm}_{1,2}$, $\lambda+ \delta_{00} = \psi_1$. Since $\xi_0^* f^*L$ defines a morphism from $\Delta_0$ to a curve, it squares to zero. On the other hand, the string equation (\cite[p.574]{ACGH2}) says that the pushforward of $\psi_1^2 = \frac{1}{c} (\xi_0^* f^*L)$ to $\bar{\shm}_{1,1}$ is the degree of $\psi = \lambda$, which we have noted is $\frac{1}{24}$, hence $\xi_0^* f^*L$ must be trivial, so that $f$ is constant. 
				\end{proof}
			
			\begin{lemma}\label{d0noa0}
				The divisor $\Delta_0$ cannot map to $\sha_0$. 
			\end{lemma}
				\begin{proof}
					We compute:
					\begin{align*}
						f^*L|_{\Delta_0} &= ((a+10b)\lambda + (c-2b)\delta_{00}\\
					\end{align*}
					Since $\Delta_0$ is contracted to a point, both the $\lambda$ and the $\delta_{00}$ factor vanish. Hence $a=-10b$ and $c=2b$, so $f^*L = -b(10\lambda - \delta_0 - 2\delta_1)$. From the relation $\delta_0 +2\delta_1 = 10\lambda$, we have $f^*L = 0$ to begin with, so that $f$ is trivial. 
				\end{proof}
			\begin{remark}
				This also holds on $\shm_{g,n}$ more generally. By \cite[Chapter 19, Lemma 4.9]{ACGH2}, the pullback map $\xi_0: \tilde{\Delta_0} \to \shm_{g,n}$ is injective on $H^2$ when $g\geq 2$. Hence $c_1(f^*L) = 0,$ so $f$ is constant. 
			\end{remark}
				
		\subsubsection{Getting Rid of Points}
			In this section, we will show that any morphism $\shm_{2,n} \to \sha_h$ for $h = 2,3$ factors through the projection to $\shm_{2}$. By the previous work, this completes the classifcation of morphisms $\shm_{2,n} \to \sha_h$ for $h =2,3$. 
			
			By \Cref{g2oneeval} we may assume that the local monodromy around the boundary divisors of $\bar{\shm}_{2,n}$ have a unique eigenvalue, either $1$ or $-1$. By twisting by a character, we may assume the local monodromy operators are unipotent. We will treat the cases $\shm_{2,n} \to \sha_2$ and $\shm_{2,n} \to \sha_3$ simultaneously. The next lemma will prove that a nonconstant morphism $\shm_{2,n} \to \sha_2$ factors throuegh the projection to $\shm_2$, while also reducing the question for a morphism $\shm_{2,n} \to \sha_3$ to the case where the local monodromy around $\Delta_0$ has nilpotent rank 2. Then, \Cref{nork2monodromy} shows that in that latter case, $f$ factors through a product locus, and is determined by a period map $\shm_{2,n} \to \sha_2$, which has already been shown to factor through the projection to $\shm_2$. 
			
			\begin{lemma}\label{nork1monodromy}
				If the local monodromy around $\Delta_0$ has nilpotent rank one, then $f$ factors through the projection to $\shm_2$.
			\end{lemma}
				\begin{proof}
					 Let $a$ and $a'$ be isotopic nonseparating simple closed curves bounding one of the marked points. The Dehn twists $T_a$ and $T_{a'}$ commute and have nilpotent rank one, hence the logarithms of their images annihilate each other. If their logarithms are scalar multiples of each other, the rank-2 braid relation calculation in \Cref{usefulcalcs} shows $f_*(T_a) = f_*(T_{a'})$, hence $f_*$ annihilates the point pushing subgroup.
					 
					  Otherwise, the images of their logarithms are different lines. Consider a nonseparating simple closed curve $b$ intersecting $a$ and $a'$ once each. Then the twists around $a$ and $a'$ satisfy the braid relation with the twist around $b$. Denote the logarithms of $f_*(T_a), f_*(T_{a'})$ and $f_*(T_b)$ by $a, a',$ and $b$, respectively. Let $v \in \img b$ and consider $w=(a-a')v$. Since $aba=-a$ and $a'ba'=-a'$, we have $bw=0$. Since $aa'=a'a=0$, $w$ is also annihilated by $a$ and $a'$. Similarly, the logarithm of every other Dehn twist generating $\Mod(2,n)$ commutes with $a$ and $a'$, hence the image of every generator of $\Mod(2,n)$ preserves $w$. By the theorem of the fixed part, the local system associated to $f$ has a constant factor. When $h=2$, this forces the map to be trivial, while when $h=3$, $f$ factors through the product locus $\sha_1 \times \sha_2$, and the map to $\sha_2$ satisfies the conditions of the previous paragraph, so that $f$ factors through the projection to $\shm_2$. 
				\end{proof}
				
			\begin{lemma}\label{nork2monodromy}
				If $f: \shm_{2,n} \to \sha_3$ is a nonconstant morphism such that the local monodromy around $\Delta_0$ has nilpotent rank two, then $f$ factors through a product locus in $\sha_3$. 
			\end{lemma}
				\begin{proof}
					We will construct an explicit invariant vector $v$ for $f_*$. If $\VV$ is the local system associated to $f$, this means $H^0(\shm_{2,n}, \VV) \neq 0$. By the theorem of the fixed part, this gives $\VV$ a constant sub-Hodge structure of dimension $2, 4$, or $6$. If the constant factor has dimension $6$, $f$ is constant. 
					
					Suppose the local monodromy around a divisor $\Delta_{1, \{p,q\}}$ is nontrivial\footnote{Recall this parameterizes a union of two genus one curves, one with $p$ marked points and one with $q$ marked points, with $p+q=n$.}. We may assume this, because if the monodromy around every such separating divisor were trivial, then $f_*$ would annihilate the Torelli subgroup and we would conclude by superrigidiy, as in the proof of \Cref{superrigidity}. Let $a$ and $b$ be nonseparating simple closed curves meeting at a point, and let $a'$ be isotopic to $a$, such that $a'$ intersects $b$ once, $a'$ and $a$ bound a marked point, and a regular neighborhood of $a' \cup b$ is a torus with $p$ marked points and one boundary component $\delta$. 
					
					Denote the logarithms of $f_*(T_x)$ by $x$. By the argument of \Cref{selfadjtrick} using symplectic adjoints, the nilpotent operator $aa'$ has rank either 2 or 0. We have that $f_*(T_\delta) -I = 12(a' + a'ba')$. Let $v \in \img(a'+a'ba')$. Then $v \in \img a'$ and $v \in \img b$, so $a'v =bv = 0$. If $\rk aa' = 2$, this also shows $av = 0$, while if $\rk aa' = 0$, then $av = 0$ to begin with. If $x$ is any nonseparating simple closed curve disjoint from $a'$, \Cref{selfadjtrick} again shows that $\rk xa'$ is $0$ or $2$. We wish to show that $xv = 0$, so suppose $\rk xa' = 2$. Then the images of $a'$ and $x$ are equal. Let $y$ be another nonseparating simple closed curve disjoint from $a'$ intersecting $x$ once. Then $T_y$ and $T_x$ satisfy the braid relation, and $T_y$ commutes with $T_{a'}$. Again write $y = \rho(T_y) -I$. Since $y$ commutes with $a'$, $y$ preserves $\img a'= \img x$. Hence $yx \subset \img x$, so $xyx=0$. By a similar argument, $yxy=0$, so $x=y$, so $\rho(T_x) = \rho(T_y)$,  and $f_*$ annihilates the commutator subgroup of $\Mod(2,n)$. We therefore must have $xa'= 0$. Since $v \in \img a'$, we conclude that $v$ is in $\ker (f_*(T_x)-I))$ for every nonseparating simple closed curve $x$ generating $\Mod(2,n)$, so $v$ is invariant under $f_*$. 
				\end{proof}
				
	\subsection{Inducting}
		We now use our classification of maps $\shm_{2,n} \to \sha_h$ as the base case for an induction on which classifies $f: \shm_{g,n} \to \sha_h$ for $h \leq 2g-1$. Suppose such maps have been classified as constant if $h<g$, and when $h \geq g$ and nonconstant, as the Torelli map, potentially factored through a product locus, for all genera $\leq g-1$ and all $n$. Our induction is on $(g,h)$ in lexicographic order, with $n$ free for each $g$.
		
		Let $f: \shm_{g,n} \to \sha_h$. By \Cref{mgagunipotent}, we may assume the local monodromy around the boundary divisors $\Delta_0$ and $\Delta_{1, \{0,n\}}$ is unipotent, so that $f$ extends to a morphism between the Deligne-Mumford and Satake compactifications. We consider the rank of the local monodromy around $\Delta_0$. Represent the local monodromy by a Dehn twist around a nonseparating simple closed curve $c$, let $\rho = f_*$, and write $\rho(T_c) = I+N$ where $N$ has order of nilpotency $2$.
		
		As we have noted before, if $N=0$, the representation is trivial, and so is the map, so we will always assume $N$ is nonzero. Moreover, we will assume $h \geq g$, as the cases $h \leq g$ have been worked out in \Cref{mgagtorelli}, and we will not rehash the argument here. Let $r = \rk N$. The monodromy weight filtration associated to $N$ is given by $0 \subset \img N \subset \ker N \subset \CC^{2h}/\ker N$. Hence the weight one graded piece $\ker N / \img N$ has rank $2h-2r$. 
		
		The limit mixed Hodge structure associated to $f$ induces a map $\phi: \shm_{g-1,n+2} \to \sha_{h-r}$.  Indeed, let $N^\times$ be a punctured tubular neighborhood of $\Delta_0$. By \cite[Proposition 2.5]{hubbardkoch}, we may identify $\pi_1(N^\times)$ with the oriented stabilizer of the isotopy class of $a$ in $\Mod(g,n)$. If $R$ the genus $g-1$ subsurface with two boundary components obtained by deleting a regular neighborhood of $c$, then $\Stab(c)^+$ is the quotient of $\Mod(R)$ obtained by identifying the Dehn twists about the two boundary components as $T_c$, and the exact sequence for $\pi_1(N^\times)$ reads
		$$1 \to \ZZ[T_c] \to \Stab(c)^+ \to \pi_1(\Delta_0) \to 1,$$ where $\pi_1(\Delta_0)$ is the quotient of $\Mod(g-1,n+2)$ by the involution exchanging the two new marked points. The restriction of $\rho$ to $\Stab(c)^+$ lies in the parabolic subgroup preserving the monodromy weight filtration, and $\pi_1(\Delta_0)$ acts on the graded pieces of the filtration. The restriction of this action to $\gr_1^W$ is the representation realizing the map $\Delta_0 \to \sha_{h-r}$. 
	 	
	 	Since $r \leq h \leq 2g-1$, the action of $\Mod(g-1,n+2)$ on $\gr_0^W$ and $\gr_2^W$ is trivial. From this, we conclude that a Dehn twist supported away from a regular neighborhood of $c$ is of the form 
		$$ \rho(T_x) = \begin{bmatrix}
			I_r & * & *\\
			0 & A & *\\
			0 & 0 & I_r\\
		\end{bmatrix}.$$
		
		We should remark that this holds even when $g=3$. While it is true that the genus 2 mapping class group has nontrivial representations of ranks $1,2,$ and $3$, those representations are all semisimple with image a quotient of $\ZZ/10\ZZ$, hence the eigenvalues of a nonseparating Dehn twist under such a representation must be roots of unity. However, we have already established that the eigenvalues of a nonseparating Dehn twist under $\rho$ are one, hence the same is true for their quotients in $\pi_1(\Delta_0)$. 
		
		Let $d$ be a nonseparating simple closed curve intersecting $c$ once. To determine the image of $\Mod(g,n)$, it is enough to determine the image of $\Stab(c)^+$ and $T_d$, since together they generate the mapping class group. 
		
		\subsubsection{\textbf{Ruling Out Higher Nilpotent Rank:}}\label{nohighernilprank} If $r \geq 2$ we have $h-r \leq 2g-3 = 2(g-1) -1$. Therefore the induced map $\phi$ on $\Delta_0$ is classified by the inductive hypothesis. In particular if $h-r<g-1$ then $\phi$ is constant, and otherwise $\phi$ factors through the projection to $\shm_{g-1},$ follwed by the Torelli map $\shm_{g-1} \to \sha_{g-1}$, and includes into $\sha_h$ by attaching a constant factor. The block $A$ in the matrix above is equal to $J(T_x) \oplus I_{2h-2r-2g}$, where $J$ is the representation corresponding to the Torelli map. We can refine the previous description and say that a Dehn twist around a nonseparating simple closed curve supported away from a neighborhood of $c$ has the form 
				$$ \rho(T_x) = \begin{bmatrix}
			I_r & * & 0 &*\\
			0 & J(T_x) & 0 & *\\
			0 & 0 & I_{2h-2g-2r} & 0\\
			0 & 0 & 0 & I_r
		\end{bmatrix}.$$
		
		There are no off-diagonal terms for $I_{2g-2g-2r}$ because the constant Abelian factor is a trivial $\Mod(g,n)$-module, and $H^1(\Mod(g,n), \CC)  =0$. 
		
		Let $\rho(T_c) - I:= c$ and $\rho(T_d) - I := d$. These are both rank $r$ nilpotent operators which square to zero. Since $T_c$ and $T_d$ satisfy the braid relation, we have $d+dcd=c+cdc$, and $\rho(T_\delta) = I+ 12(c+cdc)$. We now analyze this separating twist in detail. 
		
		Recall the rank 2 braid pair from \Cref{usefulcalcs}:
		\[U = \begin{bmatrix}
			1 & 1 \\ 0 & 1\\
		\end{bmatrix}, \quad V = \begin{bmatrix}
		1 & 0 \\ -1 & 0\\
		\end{bmatrix}.\] 
		\begin{proposition}\label{braidbasis}
			Let $\rho: \Mod(g,n) \to \GL(m, \CC)$ be a representation of the mapping class group such that the image of a Dehn twist about a nonseparating simple closed curve is unipotent, and whose logarithm $N$ is nonzero and satisfies $N^2 =0$. Let $c$ and $d$ be two nonseparating simple closed curves which intersect at a point and bound a genus one separating curve $\delta$, and suppose $\rho(T_c) \neq \rho(T_d)$. Then there exists a basis for which 
			\[
			\rho(T_c) = \begin{bmatrix}
				U^{\oplus k} & 0 & 0\\
				0 & U^{\oplus r-k} & 0\\
				0 & 0 & I
			\end{bmatrix}, \qquad 
			\rho(T_d)= \begin{bmatrix}
				U^{\oplus k} & 0 & 0\\
				0 & V^{\oplus r-k} & 0\\
				0 & 0 & I\\
			\end{bmatrix}.\]
			Here, $r = \rk (\rho(T_c) - I)$ and $k = \rk(\rho(T_\delta )-I)$.
		\end{proposition}
			\begin{proof}
				Denote by $c$ and $d$ the logarithms of $\rho(T_c)$ and $\rho(T_d)$, respectively. Since these square to zero, the braid relation implies $\rho(T_\delta) = I +12(c+cdc) = I+12(d+dcd)$. We first show $\img(c+cdc) = \img(c) \cap \img(d)$. Since $c+cdc=d+dcd$, one inclusion is clear. On the other hand if $v \in \img(c) \cap \img(d)$, write $v = cx = dy$. Then $(c+cdc)x = cx + cdcx = v + cdv = v+cd^2 y = v$, so $v \in \img(c+cdc)$. 
				
				Next, we find the subspace on which $\rho(T_c)$ and $\rho(T_d)$ both act as direct sums of the matrix $U$ and the identity. Define $H:=c-d$. We have $H^2 = -(cd+dc)$ and $H^3 = H$, so that $H$ has eigenvalues $0,1,$ and $-1$, with $H^2$ acting as a projection onto $\img(H)$. This gives  $V = \CC^{2h}$ a canonical decomposition $V = \ker(H) \oplus \img(H)$, and on $\ker(H)$, $c=d$. Hence if $v \in \ker(H),$ $cv = dv$, so $c \ker(H) \subset \img(c) \cap \img(d) = \img(c+cdc)$, while the other inclusion is immediate. 
				
				We now analyze $\img(H)$. By the braid relation, we have $(-cd)^2 = -cd$ and $(-dc)^2 = -dc$, so $-cd$ and $-dc$ are projections. Moreover, they are orthogonal, and their images span $\img(H)$, giving the decomposition $\img(H) = \img(cd) \oplus \img(dc)$. Moreover, we have decompositions $\img(c) = \img(c+dcd) \oplus \img(cd)$ and $\img(d) = \img(d+dcd) \oplus \img(dc)$. 
				
				Frp, this decomposition, we can decompose the images of $c$ and $d$ into $k$ common lines and $r-k$ lines on which $I+c$ and $I+d$ satisfy the braid relation. Indeed, $L = cdv$ for some $v$, so $dL = dcdv = -dv \subset \img d$, and since $L \not\subset \img(d+dcd)$, we in fact have $M:=dL \subset \img dc$. Once we fix a basis for $\img(cd)$, applying $d$ gives a basis for $\img(dc)$ with respect to which $\rho(T_c)$ and $\rho(T_d)$ act on $\img(H)$ as copies of the braid pair $U$ and $V$, respectively. 
			\end{proof}
		
		With this decomposition, we consider two cases: $k =0$ and $k>1$. Equivalently, these cases are when the local monodromy around the divisor $\Delta_{1,\{0,\dots, n\}}$ is either trivial or nontrivial. 
		
		\textit{(a): There is a fixed subspace if $k>0$:} 
		Using \Cref{braidbasis}, the assumption $k>0$ means there is a subspace on which $\rho(T_c)$ and $\rho(T_d)$ act nontrivially with equal images. We claim the subspace $W = \img(c+dcd)$ is preserved by $\rho(\Stab(c))$. Indeed, $\Mod(g-1,n+2)$ acts trivially on $\img(c)$ by the inductive hypothesis, and of course $\rho(T_c)$ preserves $\img(c)$. By the theorem of the fixed part, $W$ splits off as a constant isogeny factor. We conclude that $f$ factors through the inclusion $\sha_{k} \times \sha_{h-k} \to \sha_h$, and we are done by induction on $h$.

		\textit{(b): $k=0$ forces $r=1$:} We now consider when the images of $c$ and $d$ form $r$ braid pairs. When $\img c$ and $\img d$ are transverse, we know how to render $\rho(T_d)$ in the basis preserving the monodromy weight filtration induced by $c$. 
		\begin{lemma}
			Suppose $c$ and $d$ are square-zero nilpotent operators of rank $r$ with disjoint images. Relative to a basis respecting the filtration $0 \subset \img c \subset \ker c \subset \CC^{2h}$, we may write
			$$\rho(T_c) = \begin{bmatrix}
				I_r & 0 & I_r\\
				0 & I_{2h-2r} & 0\\
				0 & 0 & I_r\\
			\end{bmatrix}, \qquad \qquad \rho(T_d) = \begin{bmatrix}
			I_r & 0 & 0\\
			0 & I_{2h-2r} & 0\\
			-I_r & 0 & I_r\\
			\end{bmatrix}.$$
		\end{lemma}
			\begin{proof}
				Since $\img c \cap \img d = 0$, we have $c+cdc=d+dcd=0$, so $cdc=-c$. The decomposition established in \Cref{braidbasis} can be refined to say that $d|_{\img c}$ is an injection of one $r$ dimensional subspace into another, hence is an isomorphism. On the other hand, $\img c \simeq \CC^{2h} /\ker c$, hence $\img d \simeq \CC^{2h}/\ker c$. By symmetry, we then have $\ker c/\img c = \ker d/\img d = \ker c \cap \ker d$. Given an ingetral basis $\{v_i\}$ for $\img c$, $\{-d v_i\}$ is an integral basis for $\img d$. In this basis, $\rho(T_c) = I+c$ and $\rho(T_d) = I+d$ have the desired form. 
			\end{proof}
		Now, within $\Stab(c)^+$, consider a twist about a nonseparating curve $a$ disjoint from both $c$ and $d$. As we have established, its image is upper-triangular with respect to this basis. However, since $\rho(T_a)$ also commutes with $\rho(T_d)$, it preserves $\img d$, and there can be no off-diagonal terms. Hence 
		$$\rho(T_a) = \begin{bmatrix}
			I_r & &\\
			& J(T_a) \oplus I_{2h-2g-2r} &\\
			& & I_r
		\end{bmatrix}.$$ Since we know the image of $T_a$ has nilpotent rank one, and $T_a$ is conjugate to $T_c$, we conclude $r = 1$, a contradiction. 
		
		\subsubsection{\textbf{Transvective Monodromy:}} We now suppose $r = 1$. Then $f$ induces a map $\shm_{g-1,n+2} \to \sha_{h-1}$, which is not immediately classified by the inductive step. By the braid relation, every pair simple closed curves $(a_i,b_i)$ intersecting once act nontrivially on a plane, and trvially on the rest of $\CC^{2h}$; moreover, these $2g$ planes are disjoint, since otherwise two curves intersecting once would act nontrivially on the same line and commute. Let $a_1 =e_1, e_2 \dots, e_{n+1}$ be nonseparating simple closed curves isotopic to $a_1$ on the unmarked surface, such that each pair bounds a marked point. $\rho(e_i)$ satisfies the braid relation with $b_1$, and if the images of $\rho(T_{a_1} -I)$ and $\rho(T_{e_i}-I)$ are on the same line, then \Cref{usefulcalcs} shows that $\rho(e_i) = \rho(a_1)$, hence $\rho$ annihilates the point pushing subgroup. Similarly, if $x$ is a nonseparating simple closed curve isotopic to one of $y\in \{a_1, \dots, a_g, b_1, \dots, b_g\}$, and the images of $\rho(x)-I$ and $\rho(y)-I$ are the same line, then the braid relation implies $\rho(x) = \rho(y)$, so that the Torelli group is in the kernel of $\rho$. We conlcude in the same way as in the proof of \Cref{mgagtorelli} that $\rho$ is conjugate to the standard representation to $\Sp(2g,\ZZ)$, now with a $2h-2g$-dimensional fixed factor. 
		
		It only remains to verify the hypotheses stated in the previous paragraph: if $x$ and $y$ are isotopic simple closed curves, then the images of $\rho(T_x)-I$ and $\rho(T_y)-I$, which we will call $x$ and $y$, respectively, lie on the same line. To this end, suppose the images are different. For simplicity take $y=a_1$, so that $x$ and $a_1$ intersect $b_1$ at a point. Let $x = \rho(T_x) - I, a = \rho(T_{a_1})-I, $ and $b = \rho(T_{b_1} - I)$. Let $v \in \img b\cap \ZZ$ and consider $w=(a-x)v$. We have
		$b(a-x)v = bav - bxv = babu-bxbu = -bu +bu = 0$. Similarly $a(a-x)v = -axv$, and this vanishes since $a$ and $x$ have rank one and are assumed to have different images. By the same reasoning $x(a-x)v = 0$. Finally since $v\in \img b \subset \ker b$, $w$ is preserved by every matrix that commutes with $\rho(T_{b_1})$, which includes the images of all the other generators of $\Mod(g,n)$. Indeed, if $c$ is a nonseparating simple closed curve disjoint from $a, x,$ and $b$, write $\rho(T_c) = I+c$. Since the curves $c$ and $a$ are disjoint, Dehn twists about them commute, hence the nilpotent operators $a$ and $c$ commute. Since the images of $a$ and $c$ are disjoint lines, we must have $ac=ca=0$. It follows that $\rho(T_c)av = (I+c)av = av$. The same argument applies for $xv$, hence $(a-x)v$ is an invariant vector. By the theorem of the fixed part, $f$ has a constant sub-Hodge structure, hence the period map factors through a product locus $\sha_{h-i} \times \sha_i \to \sha_h$, and we are done by induction. 
		
		More invariantly, the argument in the previous paragraph shows that if a point-pushing map is nontrivial, its image is preserved by the entire mapping class group. After quotienting out by that image, $\rho$ factors through the projection to $\shm_g$, and for the generators of $\Mod(g)$ we have shown that there exists a basis for which $\rho$ is conjugate to the standard symplectic representation plus a trivial factor. 
		\begin{remark}
			It may be worth comparing the conclusion of our argument to that of \cite{kaufmann2026linearrepresentationsmappingclass}. After establishing unipotency and nilpotent rank $r=1$, they explicitly compute that such a representation annihilates the Johnson kernel $\mathcal{K}_g$, the kernel of the abelianization of the Torelli group. Hence $\rho(\mathcal{I}_g)$ is unipotent and abelian, so that the $\mathcal{I}_g-$fixed space of $\rho$ is nonempty. In the regime where the inductive hypothesis applies, the $\mathcal{I}_g$-fixed module must be trivial. By contrast, when $r=1$, our argument considers when a specific bounding pair map $T_aT_x^{-1}$ has nontrivial image, and shows that that image is in fact fixed by the entire mapping class group.
		\end{remark}

\printbibliography
\end{document}